%% file: 2005-40.tex
\documentclass{gtart_h}  

\input gtspec

\lognumber{583}
\received{22 March 2005}
\volumenumber{9}\papernumber{40}\volumeyear{2005}
\pagenumbers{1775}{1834}   
\revised{11 September 2005}
\published{25 September 2005}
\accepted{12 August 2005}
\proposed{Leonid Polterovich}
\seconded{Yasha Eliashberg, Eleny Ionel}

\usepackage[all]{xy}
\usepackage{amssymb,amsmath,graphicx,psfrag}

\def\psfraga <#1,#2> #3#4{%
\psfrag {#3}{\smash{\rlap{\kern #1 \raise #2\hbox{#4}}}}}

\def\figref#1{\hyperlink{#1anchor}{Figure~\ref*{#1}}}
\def\anchor#1{\noindent\hypertarget{#1anchor}{\smash{$\phantom{99}$}}}

\newcommand{\eqdef}{\;{:=}\;}

\newtheorem{Theorem}{Theorem}

\numberwithin{Theorem}{section}

\newtheorem{Lemma}[Theorem]{Lemma}
\newtheorem{Proposition}[Theorem]{Proposition}
\theoremstyle{definition}
\newtheorem{Definition}[Theorem]{Definition}
\theoremstyle{remark}
\newtheorem{Remark}[Theorem]{Remark}

\newtheorem{Corollary}[Theorem]{Corollary}

\expandafter\chardef\csname pre amssym.def at\endcsname=\the\catcode`\@
\catcode`\@=11
\def\undefine#1{\let#1\undefined}
\def\newsymbol#1#2#3#4#5{\let\next@\relax
 \ifnum#2=\@ne\let\next@\msafam@\else
 \ifnum#2=\tw@\let\next@\msbfam@\fi\fi
 \mathchardef#1="#3\next@#4#5}
\def\mathhexbox@#1#2#3{\relax
 \ifmmode\mathpalette{}{\m@th\mathchar"#1#2#3}%
 \else\leavevmode\hbox{$\m@th\mathchar"#1#2#3$}\fi}
\def\hexnumber@#1{\ifcase#1 0\or 1\or 2\or 3\or 4\or 5\or 6\or 7\or 8\or
 9\or A\or B\or C\or D\or E\or F\fi}

\font\teneufm=eufm10
\font\seveneufm=eufm7
\font\fiveeufm=eufm5
\newfam\eufmfam
\textfont\eufmfam=\teneufm
\scriptfont\eufmfam=\seveneufm
\scriptscriptfont\eufmfam=\fiveeufm

\catcode`\@=\csname pre amssym.def at\endcsname

\def    \eps    {\epsilon}

\newcommand{\Ham}{{\mathit Ham}}

\newcommand{\id}{{\mathrm{id}}}

\newcommand{\Aa}{{\mathcal A}}
\newcommand{\Bb}{{\mathcal B}}

\newcommand{\Ll}{{\mathcal L}}
\newcommand{\Jj}{{\mathcal J}}
\newcommand{\Ee}{{\mathcal E}}
\newcommand{\Ff}{{\mathcal F}}

\newcommand{\Hh}{{\mathcal H}}

\newcommand{\Mm}{{\mathcal M}}

\newcommand{\Vv}{{\mathcal V}}
\newcommand{\Ss}{{\mathcal S}}

\newcommand{\loc}{\operatorname{loc}}
\def    \wtilde  { \widetilde }
\def    \om      {{\omega}}
\def    \Om      {{\Omega}}

\def    \C      {{\mathbb C}}
\def    \R      {{\mathbb R}}
\def    \Z      {{\mathbb Z}}

\def    \T      {{\mathbb T}}

\def    \ra     {{\rightarrow}}

\def    \12    {{\textstyle\frac{1}{2}}}
\def    \HZ    {\operatorname{HZ}}
\def    \p      {\partial}
\def    \Morse    {\operatorname{\Morse}}
\def    \Maslov   {\operatorname{Maslov}}
\def    \CZ    {\operatorname{CZ}}
\def    \HZ    {\operatorname{HZ}}

\def    \Crit    {\operatorname{Crit}}
\def    \grad   {\operatorname{grad}}
\def    \test   {\operatorname{test}}
\def    \fix   {\operatorname{fix}}
\def    \codim  {\operatorname{codim}}

\def    \HF     {\operatorname{HF}}
\def    \CZ     {\operatorname{CZ}}
\def    \CF     {\operatorname{CF}}

\def    \Ham    {\operatorname{Ham}}

\newcommand{\height}{\operatorname{height}}
\newcommand{\Length}{\operatorname{Length}}
\newcommand{\Hor}{\operatorname{Hor}}
\newcommand{\volume}{\operatorname{volume}}
\newcommand{\Hopf}{\operatorname{Hopf}}
\newcommand{\xtop}{\operatorname{top}}
\newcommand{\xbottom}{\operatorname{bottom}}

\begin{document}

\title[Squeezing in Floer theory and refined HZ--capacities]{Squeezing in Floer
theory and refined\\\vglue-6pt\\Hofer--Zehnder capacities of sets near\\symplectic
submanifolds}
\covertitle{Squeezing in Floer
theory and refined\\Hofer--Zehnder capacities of 
sets near\\symplectic submanifolds}

\asciititle{Squeezing in Floer
theory and refined Hofer-Zehnder capacities of sets near symplectic
submanifolds}

\author{Ely Kerman}
\address{Mathematics, University of Illinois at Urbana--Champaign\\Urbana, IL 61801, USA\\{\rm and}\\Institute of Science, Walailak
University\\Nakhon Si Thammarat, 80160, Thailand}
\asciiaddress{Mathematics, University of Illinois at Urbana-Champaign\\Urbana, IL 61801, USA\\and\\Institute of Science, Walailak
University\\Nakhon Si Thammarat, 80160, Thailand}
\email{ekerman@math.uiuc.edu}
\urladdr{http://www.math.uiuc.edu/~ekerman/}
\asciiurl{http://www.math.uiuc.edu/ ekerman/}

\begin{abstract}
We use Floer homology to study the Hofer--Zehnder
capacity of neighborhoods near a closed symplectic submanifold $M$ of a
geometrically bounded and symplectically aspherical ambient manifold.
We prove that, when the unit normal bundle of $M$ is homologically
trivial in degree $\dim(M)$ (for example, if $\codim(M) > \dim(M)$),
a refined version of the Hofer--Zehnder capacity is
finite for all open sets close enough to $M$.
We compute this capacity for certain tubular
neighborhoods of $M$ by using a squeezing argument in
which the algebraic framework of Floer theory
is used to detect nontrivial periodic orbits \cite{fh,fhw,gg}.
As an application, we partially recover some existence 
results of Arnold \cite{ar} for
Hamiltonian flows which describe a charged particle moving in a
nondegenerate magnetic field on a torus. Following \cite{kl}, we also
relate our refined capacity to the study of Hamiltonian paths with
minimal Hofer length.
\end{abstract}

\asciiabstract{%
We use Floer homology to study the Hofer-Zehnder capacity of
neighborhoods near a closed symplectic submanifold M of a
geometrically bounded and symplectically aspherical ambient manifold.
We prove that, when the unit normal bundle of M is homologically
trivial in degree dim(M) (for example, if codim(M) > dim(M)), a
refined version of the Hofer-Zehnder capacity is finite for all open
sets close enough to M.  We compute this capacity for certain tubular
neighborhoods of M by using a squeezing argument in which the
algebraic framework of Floer theory is used to detect nontrivial
periodic orbits.  As an application, we partially recover some
existence results of Arnold for Hamiltonian flows which describe a
charged particle moving in a nondegenerate magnetic field on a torus.
We also relate our refined capacity to the study of Hamiltonian paths
with minimal Hofer length.}


\asciikeywords{Hofer-Zehnder capacity, symplectic submanifold, Floer
homology}

\maketitle

\section{Introduction and results}

The Hofer--Zehnder capacity of an open subset of a symplectic
manifold, is a measure of the size of the set in terms of the
periodic orbits of the autonomous Hamiltonian flows it supports.
This symplectic invariant was introduced by Hofer and Zehnder in
\cite{hz:cap}, where the authors also compute its value for the
standard symplectic cylinder $Z^{2l}(R) = B^2(R) \times \R^{2l-2}$
in $(\R^{2l},\Omega_{2l})$. In particular, they prove that
\begin{equation}
\label{cylinder}
 c_{\HZ}(Z^{2l}(R)) = \pi R^2.
\end{equation}
This remarkable result illustrates the symplectic nature of the
Hofer--Zehnder capacity by distinguishing it from the volume.
Moreover, it yields an alternative proof of Gromov's famous
nonsqueezing theorem and hence it serves as an important link
between the existence problem for periodic orbits of Hamiltonian
flows and symplectic rigidity phenomena \cite{hz:book}.

Viewing the symplectic cylinder $Z^{2l}(R)$ as a tubular
neighborhood of the symplectic submanifold $\R^{2l-2} \subset
\R^{2l}$, equation \eqref{cylinder} naturally leads one to the
following question which motivates this work.

\textsl{Let M be a closed symplectic
submanifold of a symplectic manifold $(W,\Om)$, and let $U_R$ be a
symplectic tubular neighborhood of $M$ with (sufficiently small)
radius $R$. Is the Hofer--Zehnder capacity of $U_R$ equal to $\pi
R^2$?}

Of course, to make this precise one must specify what
is meant by a symplectic tubular neighborhood of radius $R$. These
are neighborhoods of $M$ on which the symplectic form has a standard
normal form. They are defined below in Section~\ref{sec:U_R}.

It follows from Weinstein's Symplectic Neighborhood Theorem that any
symplectic invariant of $U_R$ is completely determined by three
factors; the radius $R$, the restriction of $\Om$ to $M$, and the
isomorphism class of the normal bundle of $M$ in $W$. The real
question then is whether the Hofer--Zehnder capacity is independent
of the last two factors. It is interesting to note that an
affirmative answer in either case would further distinguish the
Hofer--Zehnder capacity from the volume. In particular, Hermann
Weyl's famous formulas for the volumes of {\em tubes}
\cite{gr,weyl}, demonstrate that the volume of a tubular
neighborhood of a submanifold depends, in general, on both the
volume of the submanifold as well as the class of its normal bundle.

The question above is also relevant to several active areas of
research in Hamiltonian dynamics and symplectic topology. In the
study of the Hamiltonian flows which describe the classical motion
of a charged particle in an electro-magnetic field, an affirmative
answer to this question immediately implies the existence of
periodic orbits on almost all low energy levels (see, for example,
\cite{cgk,gg,schl} and Section~\ref{classical} below). Following the work
of Lalonde and McDuff from \cite{lm}, this question is also of great
importance in the study of Hamiltonian paths which minimize the
Hofer length \cite{mcsl}. As well, Biran's decomposition theorem
from \cite{bi} shows that a closed K\"{a}hler manifolds can, in some
suitable sense, be ``filled-up" by a symplectic tubular
neighborhood. Hence, an affirmative answer to the question above
would have many implications for the symplectic topology of closed
K\"{a}hler manifolds.

The first results concerning this question were obtained for trivial
symplectic tubular neighborhoods of the form $M \times B^2(R)$,
where $M$ is closed and the product is given the obvious split
symplectic form. These began with the following theorem.
\begin{Theorem}[Floer--Hofer--Viterbo \cite{fhv}] If $(M, \om)$ is
weakly exact, that is, if $\omega$ vanishes on $\pi_2(M)$, then $c_{\HZ}(M
\times B^2(R)) = \pi R^2$.
\end{Theorem}
\noindent This was later extended to the case when $(M, \omega)$ is
rational by Hofer and Viterbo in \cite{hv}. In the context of
studying length minimizing Hamiltonian paths, McDuff and Slimowitz
also proved more general results for slightly different capacities
in \cite{mcsl}. The result for general symplectic manifolds $(M,
\omega)$ was proved recently by G. Lu in \cite{lu2} using the theory
of Gromov--Witten invariants and, in particular,
Liu and Tian's construction of an equivariant virtual moduli cycle from
\cite{lt}.


In this paper, we prove that a refined Hofer--Zehnder capacity of a
symplectic tubular neighborhood $U_R$ is equal to $\pi R^2$ for all
small $R>0$, provided that the ambient manifold $(W, \Om)$ is
geometrically bounded and symplectically aspherical, and the
homology of the unit normal bundle of $M$ splits in degree $2m =
\dim(M)$. To state the result precisely, we must first recall the
definition of the Hofer--Zehnder capacity and introduce the
refinements considered here.

\begin{Remark}
Several other recent works \cite{bi,cgk,gg,ma,schl} have also been
devoted to this question. A comparison of our results and methods to
those from \cite{bi,cgk,gg,ma,schl} is contained in
Sections~\ref{sec:compare} and \ref{sec:methods} below.
\end{Remark}

\subsection{Hofer--Zehnder capacities}

 On a symplectic manifold $(W, \Omega)$ every
Hamiltonian $H \co W \to \R$ determines a unique Hamiltonian vector
field $X_H$ via the equation
$$
\Om(X_H, \cdot) = -dH(\cdot).
$$
If $dH$ is compactly supported, then the flow of $X_H$ is defined
for all time (we will also refer to this as the Hamiltonian flow of $H$).
We denote the flow for time $t$ by $\phi^t_H$.

Let $U$ be a nonempty open subset of $W$. The set of \emph{test
Hamiltonians} on $U$, $\Hh_{\test}(U)$, is defined to be the set of
smooth nonnegative functions $H \co W \to [0,\,\infty)$ with the
following properties:
\begin{itemize}
    \item $H$ has compact support contained in $U$.
    \item The set on which $H$ takes the value zero has
nonempty interior.
    \item The set on which $H$ attains its maximum value,
    $\max(H)$, has nonempty interior.
\end{itemize}
A test Hamiltonian $H$ is said to be \emph{admissible} if:
\begin{description}
    \item[$(\ast)$] \textsl{The Hamiltonian flow of $H$ has no nonconstant
    periodic orbits with period less than or equal to one.}
\end{description}
The Hofer--Zehnder capacity of $U$ in $(W, \Om)$ is then defined to
be
\begin{equation*}
    c_{\HZ}(U) = \sup \left\{\max(H) \mid H \in \Hh_{\test}(U),\,\,H
    \text{  satisfies  }(\ast)\right\}.
\end{equation*}
By changing the admissibility criterion and/or the set of test
Hamiltonians one obtains different Hofer--Zehnder capacities,
\cite{gg,lu,ma,sc}. Here we consider a weaker admissibility
criterion which is defined using the action functional on $\Ll$, the
space of smooth contractible loops in $W$. For this reason, we will
assume that $(W, \Om)$ is weakly exact, ie $\Om|_{\pi_2(W)} =0.$
Under this assumption, each smooth map $x \co S^1 \to W$ has a
well-defined \emph{symplectic area}
$$
\Aa(x)= \int_{D^2} v^* \omega,
$$
where $v$ is any smooth map from the disc $D^2$ to $W$ which
satisfies $v|_{\partial D^2}=x$. For a fixed $H \in \Hh_{\test}(U)$,
the \emph{action} of $x \in \Ll$ is then defined by
$$
\Aa_H(x)= \int_0^1 H(x(t)) \, dt - \Aa(x).
$$
The critical points of the action functional $\Aa_H \co \Ll \to \R$
are precisely the contractible periodic orbits of $H$ with period
equal to one. A simple, but crucial, observation is that the
$\Aa_H$--value of any constant periodic orbit of $H$ must be less
than or equal to $\max(H)$. This leads us to the following
admissibility criteria for Hamiltonians on weakly exact manifolds.
\begin{description}
    \item[$(\ast)^{1,\,\kappa}$] \textsl{The Hamiltonian flow of $H$ has no
    nonconstant contractible periodic orbits with period equal to one and action
    in the interval
$$(\max(H), \max(H) +\kappa)$$ for $\kappa>0$.}
\end{description}

The corresponding capacity of an open subset $U \subset W$ is then
defined by
\begin{equation*}
    c_{\HZ}^{\,1,\,\kappa}(U)= \sup \left\{\max(H) \mid H \in
    \Hh_{\test}(U), H \text{  satisfies  }(\ast)^{1,\,\kappa}\right\}.
\end{equation*}
This refined capacity contains the maximum amount of information
yielded by the squeezing techniques we employ here. The effort
invested in retaining this information is repaid by the fact that
the almost existence theorem for this capacity includes a bound for
the symplectic area of the periodic orbits, (Theorem
\ref{thm:ae-area}). Similar area bounds were obtained by Hofer and
Zehnder in their original work \cite{hz:wc} on dense existence results in
$(\R^{2l}, \Om_{2l})$. These bounds, as well as some
applications, are described in Section~\ref{sec:app}.

The invariant $c_{\HZ}^{\,1,\,\kappa}$ has the usual properties of a
symplectic capacity:
\begin{enumerate}
\item {\bf Invariance}\qua If $\phi$ is a symplectomorphism of $(W,\Om)$,
then $$c_{\HZ}^{\,1,\,\kappa}(U) = c_{\HZ}^{\,1,\,\kappa}(\phi(U)).$$
\item {\bf Monotonicity}\qua If $U \subset V$, then $c_{\HZ}^{\,1,\,\kappa}
(U) \leq c_{\HZ}^{\,1,\,\kappa}(V)$.
\item {\bf Homogeneity}\qua Let the notation $c_{\HZ}^{\,1,\,\kappa}(U,\Omega)$
reflect the dependence of the capacity on the symplectic form. Then
for any constant $\alpha > 0$
$$c_{\HZ}^{\,1,\,\kappa}(U,\alpha \Omega) = |\alpha|c_{\HZ}^{\,1,\,
\kappa/|\alpha|}(U,\Omega).$$
\item {\bf Normalization}\qua If $(W, \Om)= (\R^{2l}, \Om_{2l})$, then
$$c_{\HZ}^{\,1,\,2\pi R^2}(B^{2l}(R)) = \pi R^2 = c_{\HZ}^{\,1,\,2\pi R^2}(Z^{2l}(R)).$$
\end{enumerate}
The first three of these properties follow easily from the
definition \cite{hz:book}. The normalization condition, in
particular the fact that $$c_{\HZ}^{\,1,\,2\pi R^2}(Z^{2l}(R))=\pi
R^2,$$ can be proved using the methods developed here.

We also note that for all $\kappa'> \kappa>0 $ we have
$$
c_{\HZ}^{\,1,\,\kappa}(U) \geq c_{\HZ}^{\,1,\,\kappa'} \geq
c_{\HZ}(U).
$$
Moreover, since our capacity only detects contractible periodic
orbits, there are sets $U$ with $ c_{\HZ}^{\,1,\,\kappa}(U)= \infty$
and $c_{\HZ}(U)< \infty$. (For example, let $U$ be a neighborhood of
a noncontractible loop in $\T^2$.)

\subsection{The main result}

For the remainder of the paper our ambient symplectic manifold will
be $(W, \Omega)$, and $M \subset W$ will be a closed symplectic
submanifold of dimension $2m$ and codimension $2n$. The restriction
of $\Omega$ to $M$ will be denoted by $\omega$.

We will assume that $(W, \Om)$ is \emph{symplectically aspherical}.
That is, $\Om|_{\pi_2(W)}=0$ and $c_1|_{\pi_2(W)}=0$, where the
forms $\Omega$ and $c_1$ act on $\pi_2(W)$ by integration over
(piecewise) smooth representatives.

We also assume that $(W, \Om)$  is \emph{geometrically bounded},
which means that $W$ admits an almost complex structure $J$ and a
complete Riemannian metric $g$ for which
\begin{itemize}
\item there are positive constants $c_1$ and $c_2$ such that
$$
\Om(X, JX) \geq c_1 \|X\|^2 \,\,\text{  and   }\,\, |\Om(X,Y)| \leq
c_2 \|X\|\,\|Y\|
$$
for all $X,Y \in TW$,
\item the sectional curvature of $(W,g)$ is bounded from above and the injectivity radius
of $(W,g)$ is bounded away from zero.
\end{itemize}
This assumption is trivial when $W$ is compact. For noncompact
manifolds it implies the required $C^0$--bounds for the Floer moduli
spaces.

Finally, let $E \to M$ be a normal bundle of $M$ in $W$, and let
$S(E)$ be the corresponding unit normal bundle with respect to some
fibrewise metric on $E$. We say that $S(E)$ is \emph{homologically
trivial in degree $2m$}, if the homology of $S(E)$ (with
coefficients in $\Z_2$) splits in degree $2m$, ie
$$
H_{2m}(S(E), \Z_2)= H_{2m}(M, \Z_2) \oplus
 H_{2(m-n)+1}(M, \Z_2).\footnote{
For simplicity, group isomorphisms  will always be denoted with an
``$=$" sign.}
$$
For example, this splitting occurs if the normal bundle of $M$
admits a nonvanishing section. It follows from the Gysin sequence
that this condition is satisfied if the codimension of $M \subset W$
is greater than its dimension $(n>m)$, (see, for example,
\cite[Section~11]{bt}).

The following is our main result:

\begin{Theorem}
\label{thm:main} Suppose that $(W, \Omega)$ is geometrically
bounded and symplectically aspherical. Let $M \subset W$ be a closed
symplectic submanifold whose unit normal bundle is homologically
trivial in degree $2m$. Then, for all sufficiently small $R>0$, the
symplectic tubular neighborhood $U_R$ has refined capacity
$$c_{\HZ}^{\,1,\,2 \pi
R^2}(U_R) = \pi R^2.$$
\end{Theorem}

The monotonicity property of the capacity implies the following corollary.

\begin{Corollary}\label{cor2}
For $(W, \Om)$ and $M$ as in Theorem \ref{thm:main}, every open set
$U \subset U_R \subset W$ has $c_{\HZ}^{\,1,\,2\pi R^2}(U)\leq \pi
R^2$.
\end{Corollary}

\begin{Remark}
The size of the values of $R$ for which our method of proof for Theorem
\ref{thm:main} works is restricted by two factors. First, we are
required to consider small enough neighborhoods of $M$ on which
$\Omega$ has a standard normal form. More importantly, even for
trivial tubular neighborhoods, we need $\pi R^2$ to be less than the
Gromov width of $(M,\omega)$.
\end{Remark}

\subsection{Other recent results} \label{sec:compare} We now describe
the content of some other recent works \cite{bi,cgk,gg,ma,schl}
which also consider the capacity of symplectic tubular
neighborhoods.

As part of his work in \cite{bi}, Biran finds an upper bound for the
Gromov width of symplectic tubular neighborhoods of certain closed
codimension-two symplectic submanifolds of closed K\"{a}hler
manifolds.

The papers \cite{cgk,gg,ma} all consider various {\emph{relative}}
capacities of symplectic tubular neighborhoods. The monotonicity
property of these relative capacities applies to small open sets
which \emph{contain} $M$. This is adequate to obtain new
almost/dense existence results for certain Hamiltonian flows which
describe the classical motion of a charged particle in a
nondegenerate magnetic field (see Section~\ref{classical}).

In \cite{cgk} and \cite{gg}, the ambient manifold $(W, \Om)$ is
assumed to be geometrically bounded and symplectically aspherical.
The paper \cite{cgk} considers a relative homological symplectic
capacity which is defined as in \cite{bps}. For symplectic tubular
neighborhoods $U_R$ with sufficiently small radius, it is shown that
this relative capacity equals $\pi R^2$. In \cite{gg}, the authors
consider a relative version of the Hofer--Zehnder capacity which is
defined using test Hamiltonians that attain their maximum value in
an open neighborhood of the symplectic submanifold. They prove that
the relative Hofer--Zehnder capacity of $U_R$ is also equal to $\pi
R^2$. This allows the authors to improve the dense existence results
from \cite{cgk} to almost existence results.

In \cite{ma}, Macarini introduces a stabilization procedure which
allows him to bound the capacity of a symplectic manifold $(W, \Om)$
which admits a free Hamiltonian circle action in terms of the
capacity of compact subsets of  $W \times T^*S^1$. These sets
can be symplectically embedded into a trivial symplectic tubular
neighborhood $W \times B^2(R)$ whose capacity is known to be finite
in various cases by \cite{hv,lu,lu2,mcsl}. When applied to $U_R \setminus
M$ this yields bounds on certain relative capacities and allows Macarini
to further improve the existence results from \cite{cgk,gg} by relaxing the
assumption that the ambient manifold is symplectically aspherical.

Most recently, in \cite{schl}, Schlenk has obtained a powerful new
generalization of Hofer's energy--capacity inequality from
\cite{ho2}, where the standard capacity is replaced by the
contractible capacity and the displacement energy is replaced by the
stable displacement energy. The proof involves Macarini's
stabilization procedure as well as several techniques and theorems
from the study of Hofer's geometry on the space of Hamiltonian
diffeomorphisms.

Schlenk uses his energy--capacity inequality, together with the
results of Laudenbach, Polterovich and Sikorav on the displacement
of subsets of symplectic manifolds, to obtain the following result:
\begin{Theorem}[Schlenk \cite{schl}]
\label{thm:schlenk} Let $(W, \Om)$ be a symplectic manifold which is
geometrically bounded and stably strongly semi-positive. Let $M$ be
any closed submanifold of $W$ such that either $\dim M < \codim M$,
or $\dim M = \codim M$ and $M$ is not Lagrangian. Then every
sufficiently small tubular neighborhood of $M$ has finite
(contractible) Hofer--Zehnder capacity.
\end{Theorem}

The hypotheses of this theorem are surprisingly mild and
consequently, when $\dim M \leq \codim M$, it has much broader
applications to Hamiltonian dynamics than one obtains from the
bounds on capacities from  \cite{cgk,gg,ma} or Theorem
\ref{thm:main}.

On the other hand, the more restrictive assumptions of Theorem
\ref{thm:main} allow us to use the refined capacity
$c_{\HZ}^{1,\kappa}$. The almost existence theorem for this capacity
includes bounds on the symplectic area of periodic orbits (Theorem
\ref{thm:ae-area}), and in certain cases these bounds can be used
to obtain stronger existence theorems, (Lemma \ref{lem:apriori}).

\subsection{Comparison of techniques}
\label{sec:methods}

Outside of the techniques developed by Hofer and Zehnder for
$(\R^{2l}, \Omega_{2l})$, there are three general methods known to
the author for bounding Hofer--Zehnder capacities.

The first of these methods was introduced in \cite{fhv} and
developed in \cite{hv} (see also \cite{lt,lu,lu2,mcsl}). This method
utilizes Floer's version of Gromov's compactness theorem for
perturbed $J$--holomorphic spheres. Roughly speaking, one starts
with a test Hamiltonian $H$ and a nonempty moduli space $\Ss$ of
regular $J$--holomorphic spheres with marked points which get mapped
to $x_0 \in \{H=0\}$ and $x_{\max} \in \{H=\max(H)\}$. Then, the
test Hamiltonian $H$ is used to produce the family of perturbed
Cauchy--Riemann equations $\overline{\partial}_J + \nabla(\lambda H)=0$,
where $\lambda \geq 0$. One can then show that: the collection of
moduli spaces $(\lambda,\Ss_{\lambda})$ is noncompact;
$\Ss_{\lambda}$ is nonempty for small $\lambda>0$; and
$\Ss_{\lambda}$ is empty for all $\lambda$ greater than some
positive constant (which bounds the capacity from above). Under
suitable assumptions the source of noncompactness can then be
identified as convergence to a broken Floer cylinder which breaks
along a nonconstant orbit of some $\lambda H$.

In the present work, we use a squeezing argument in Floer homology
from \cite{gg}. This is complementary to the previous method since
it works best in the absence of holomorphic spheres and relies on
Floer's gluing theorem as well as his compactness result. The
general idea is to use the algebraic framework of filtered Floer
homology to detect nontrivial orbits. The origins of the strategy
can be found in the computations of Symplectic Homology in
\cite{fhw} (see also \cite{bps,cgk}). In \cite{gg}, these ideas were
distilled by Ginzburg and G\"{u}rel into an effective method for
finding upper bounds for Hofer--Zehnder capacities.

To illustrate the basic principle, consider a test function $H$ on a
set $U$. The function
 $H$ can be approximated from above and below by two model Hamiltonians
$G_+$ and $G_-$, respectively, whose dynamics (and Floer homology)
is completely understood. One then considers the monotone Floer
continuation maps
$$
\sigma_{G_- G_+} \co \HF^{a,\,b}(G_+) \to \HF^{a,\,b}(G_-),
$$
for action values $a > \max(H)$. Floer's compactness and gluing
theorems imply that this map must factor through $\HF^{a,\,b}(H)$.
Hence, if one can show that $\sigma_{G_- G_+}$ is nontrivial, then
it follows that $\HF^{a,\,b}(H)\neq 0 $ and $H$ must have a
$1$--periodic orbit with action greater than $\max(H)$. As described
above, such an orbit must be nonconstant. In this way, one obtains
an algebraic version of the existence problem. However, there are
simple examples where the map $\sigma_{G_- G_+}$ is trivial for all
$a > \max(H)$. Moreover, even when $\sigma_{G_- G_+}$ is nontrivial,
the precise transfer mechanism is often quite subtle.

The results from \cite{bps,cgk,fhw,gg}, are all established using
this type of argument. In these papers, it is possible for the
authors to work at the homology level because the model Hamiltonians
$G_+$ and $G_-$ have qualitatively similar dynamics on the relevant
level sets. In the present work, this is not the case (the relevant
level sets are not even diffeomorphic). To overcome this difficulty
we are required to also work at the chain level in Floer theory.

For a description of the methods from Hofer's geometry used by
Schlenk in \cite{schl}, the reader is referred to the paper itself
as well as to the references \cite{hz:book, po:book}. We only
mention here that at the center of the argument lies the main result
from \cite{mcsl} which gives a sufficient condition, in terms of
periodic orbits, for a Hamiltonian $H$ to generate a length
minimizing path with respect to Hofer's norm. In turn, the main
result in \cite{mcsl} is implied by  an upper bound for a particular
Hofer--Zehnder capacity of {\em trivial} symplectic tubular
neighborhoods \cite{lm}. It is a remarkable fact that all the
finiteness theorems for capacities implied by Schlenk's Theorem
\ref{thm:schlenk} can be traced to this one bound for the capacity
of trivial symplectic tubular neighborhoods.

\subsection{Extensions of Theorem~\ref{thm:main}}

We now describe how some of the hypotheses of Theorem \ref{thm:main}
can be relaxed.

Most importantly, it should be possible to relax the hypothesis that
$(W, \Om)$ is symplectically aspherical and instead assume that the
minimal Chern number is sufficiently large and there exists a
symplectic area gap, ie
\begin{equation*}
\label{gap}
 \inf_{A \in \pi_2(W)}\biggl\{\Bigl|\int_A \Om\Bigr| : \int_A \Om
\neq 0\biggr\}
> 0.
\end{equation*}
One can then repeat the chain level arguments described here using
the Floer--Novikov complex from \cite{hs}, which also comes with a
filtration by the action functional. In particular, as described
above, the proof of Theorem \ref{thm:main} relies on an analysis of
the Floer complexes of two model Hamiltonians restricted to a
specific action interval, say $(a,\,b)$. The orbits of these model
Hamiltonians come with natural choices of spanning discs. By
choosing $b-a$ to be sufficiently small, the area gap (and the
appropriate lower bound for the minimal Chern number) allows one to
ignore all other choices of spanning discs.

It will also be clear from the proof that Theorem \ref{thm:main}
holds for the Hofer--Zehnder capacity which is defined using time
dependent test Hamiltonians $H \co [0,1] \times W \to \R$ such that
$H$ has compact support in $[0,1]\times U_R$ and each function $H_t
= H(t,\cdot)$ attains its maximum value on a common set with
nonempty interior. (Such functions play a prominent role in
\cite{sc}).


\begin{Remark}
The possible applications of the squeezing method for computing new
bounds for Hofer--Zehnder capacities are restricted primarily by the
fact that standard Floer theory counts only perturbed holomorphic
cylinders (twice punctured spheres) and hence can only see periodic
orbits in one homotopy class. In a joint project with V\,L Ginzburg
and B G\"{u}rel \cite{ggk} we will construct a more versatile
version of Floer theory following the blueprint of Symplectic Field
Theory \cite{egh}.
\end{Remark}

\subsection{Organization of the paper}

In the next section, we describe some applications of Theorem
\ref{thm:main} to classical Hamiltonian flows and Hofer's geometry.
The remainder of the paper is devoted to the proof of Theorem
\ref{thm:main}. In Section~3, we define symplectic tubular
neighborhoods and reduce Theorem \ref{thm:main} to a result about
the dynamics of test Hamiltonians on these neighborhoods, Theorem
\ref{thm:main-1dyn}. The necessary tools from Floer theory are then
described in Section~4. In Section~5, we prove Theorem \ref{thm:main-1dyn}
by reducing it to a sequence of statements about the nontriviality
of monotone Floer continuation maps which are then proved using the
tools from Section~4. Appendix~\ref{app:a} contains the proof of an upper bound
for the area of planar curves with positive curvature. This bound is
needed for the applications in Section~2. Finally, a proof of a
well-known result for the Morse homology of fibre-bundles is
contained in Appendix~\ref{app:b}.

\subsection{Acknowledgements}

The author would like to express his profound gratitude to Viktor
Ginzburg for many valuable conversations concerning this work. He is
also grateful to Stephanie Alexander, Dusa McDuff, and Felix Schlenk
for their helpful comments and suggestions. The author would also
like to thank the referee of this paper for the extremely thorough
and insightful report. Most of this work was undertaken and
completed during the author's stays at Stony Brook University and
Walailak University. He thanks these institutions for their warm
hospitality and generous support. This research was partially
supported by NSF Grant DMS--0405994.

\section{Applications}
\label{sec:app}

\subsection{An almost existence theorem with area bounds}
\label{sec:ae}

When one can prove that the Hofer--Zehnder capacity of a symplectic
manifold is finite, there are remarkable consequences for the
autonomous Hamiltonian flows on the manifold. In particular, one
obtains the {\em almost existence theorem} of Hofer--Zehnder and
Struwe \cite{hz:cap,st}. We recall the statement of this result
below and show that the corresponding result for the refined
capacity $c_{\HZ}^{1,\kappa}$ includes bounds on the symplectic area
of periodic orbits. In fact, we will use a slight generalization of
the almost existence theorem due to Macarini and Schlenk
\cite{masc}.

Let $S \subset (W,\Om)$ be a hypersurface, ie an oriented closed
submanifold of codimension one. The restriction of $\Om$ to $S$ has
a one-dimensional kernel which determines the  \emph{characteristic
foliation} of $S$. A \emph{closed characteristic} of $S$ is an
embedding of the circle into $S$ whose image is a closed leaf of
this foliation.

If $S$ is a regular level set of a Hamiltonian $H$, then the
periodic orbits of $H$ on $S$ are in one-to-one correspondence with
the closed characteristics of $S$. As demonstrated by the
counterexamples in \cite{gi-ce1,gi-ce2,gg-ce,he,ke-ce}, not every
regular level set (hypersurface) carries a periodic orbit (closed
characteristic). However it is still fruitful to study the existence
question {\em near} a fixed hypersurface $S$. Following
\cite{hz:wc}, one does this by considering a \emph{thickening} of
$S$. This is a diffeomorphism $\psi \co (-1,1) \times S \to W$ onto
an open and bounded neighborhood of $S$ such that $\psi(0, \cdot)
\co S \to W$ is the inclusion map. In other words, $\psi$ determines
a family of hypersurfaces $S_r = \psi(\{r\} \times S)$ which are
modeled on $S$.

\begin{Theorem}[Hofer--Zehnder \cite{hz:cap}, Macarini--Schlenk
\cite{masc} and Struwe \cite{st}]
Let $U$ be a neighborhood of a hypersurface
$S$ in $(W, \Om)$ such that $c_{\HZ}(U) < \infty$. If  $\psi$ is a
thickening of $S$ whose image is contained in $U$, then there is a
closed characteristic on $S_r$ for almost every $r \in (-1,1)$.
\end{Theorem}

The almost existence theorem corresponding to the refined capacity
$c_{\HZ}^{1,\kappa}$ includes bounds on the symplectic area.
\begin{Theorem}
\label{thm:ae-area} Let $S$ be a hypersurface of a weakly exact
symplectic manifold $(W, \Om)$, and let $U$ be a neighborhood of $S$
such that  $c_{\HZ}^{1,\kappa}(U) < \infty$ for some $\kappa <
\infty$. If $\psi$ is a thickening of $S$ with image in $U$, then
almost every hypersurface $S_r$ carries a contractible closed
characteristic $x$ with symplectic area satisfying  $|\Aa(x)| \leq
c_{\HZ}^{1,\kappa}(U) +\kappa$.
\end{Theorem}

\proof
This result follows easily from the arguments in \cite{hz:cap,st}
and their refinements in \cite{masc}. The key point is that during
the usual limit process, outlined below,  one can use the bounds on
the action given by the capacity $c_{\HZ}^{1,\kappa}$ to obtain the
bounds on the symplectic area.

The limit process from \cite{hz:cap,masc,st} is structured as
follows. For a fixed $r \in (-1,1)$ one constructs a sequence of
test Hamiltonians $F_i$ on $U$ such that $\max(F_i)
> c_{\HZ}(U)$. It follows that each $F_i$ has a nonconstant periodic
orbit $y_i$ on some level $S_{r_i}$. The functions $F_i$ are
constructed so that the $r_i$ converge to $r$. To prove the almost
existence theorem, one then shows that the $y_i$ converge to a
periodic orbit $x$ on $S_r$ for almost all values of $r$.

If this procedure is repeated using the refined capacity
$c_{\HZ}^{1,\kappa}$, we obtain a sequence of nonconstant
contractible periodic orbits $y_i$ with period equal to one, such
that
$$
\max(F_i) < \int_0^1 F_i(y_i(t))\,dt -\Aa(y_i) < \max(F_i) + \kappa.
$$
Since each $F_i$ is nonnegative and $\int_0^1 F_i(y_i(t))\,dt <
\max(F_i)$, this implies that
$$
0 > \Aa(y_i) > - \max(F_i) - \kappa.
$$
The $F_i$ can also be chosen so that $\max(F_i) -
c_{\HZ}^{1,\kappa}>0$ is arbitrarily small. In particular, for a
sequence $\epsilon_i \to 0^+$ we have
$$
0 > \Aa(y_i) \geq - c_{\HZ}^{1,\kappa}- \kappa - \epsilon_i.
$$
Now, when the sequence of orbits $y_i$ converges, the limit $x
\subset S_r$ is a nontrivial contractible periodic orbit with period
equal to one and symplectic area satisfying
$$
0 \geq \Aa(x) \geq - c_{\HZ}^{1,\kappa}- \kappa.\eqno{\qed}
$$

The existence question for thickenings of a hypersurface was first
considered for $(\R^{2l}, \Om_{2l})$ by Hofer and Zehnder in
\cite{hz:wc}. There the authors establish the existence of closed
characteristics on a dense set of hypersurfaces in a thickening, and
they also obtain the same information about the symplectic area of
these characteristics. Since \cite{hz:wc} predates Floer's creation
of his homology theory, it goes without saying that the methods used
in \cite{hz:wc} are much different from those used here.

As noted in \cite{hz:wc}, the bounds on the symplectic area can be
used, in the presence of certain \emph{a priori} bounds, to prove
the existence of periodic orbits on {\em fixed} hypersurfaces (see
also \cite{br,bhr}). More precisely, consider a situation where
Theorem \ref{thm:ae-area} holds. That is, we have a hypersurface $S$
and a constant $K>0$ such that for any suitable thickening $\psi$
almost every hypersurface $S_r$ of $\psi$ contains a contractible
closed characteristic $x$ with $|\Aa(x)| \leq K$.

\begin{Lemma}{\rm\cite{hz:wc}}\qua
\label{lem:apriori}
Suppose there is a thickening $\psi$ of $S$
and a constant $C>0$ such that for any closed characteristic $y$
in the image of $\psi$ we have
$$
|\Aa(y)| \geq Cl(y),
$$
where $l(y)$ is the length of $y$ with respect to some fixed
metric. Then $S$ carries a contractible closed characteristic $x$
 with $|\Aa(x)| \leq  c_{\HZ}^{1,\kappa}(U) +\kappa$.
\end{Lemma}

\begin{proof}
By Theorem \ref{thm:ae-area} we can find a sequence $r_i \to 0$ such
that there is a closed characteristic $y_i$ on $S_{r_i}$ whose areas
are uniformly bounded from above. The {\em a priori} bounds in Lemma
\ref{lem:apriori} yield a uniform upper bound for the lengths of the
$y_i$. Since the lengths of closed characteristics are also
uniformly bounded away from zero \cite{hz:wc}, the result follows
immediately from the Arzela--Ascoli Theorem.
\end{proof}

If $S$ is a regular level set of some Hamiltonian $H$ and the
thickening $\psi$ is determined by the nearby level sets of $H$,
then the bound in Lemma \ref{lem:apriori} can be replaced by one
involving the period instead of the length. In particular it
suffices to find a constant $C$ such that for any periodic orbit $y$
of $H$ with period $T$ we have
$$
|\Aa(y)| \geq CT.
$$

\subsection{Classical Hamiltonian flows}
\label{classical}

Using Theorem \ref{thm:ae-area} and Lemma \ref{lem:apriori}, we
describe in this section  some applications of Theorem
\ref{thm:main} to classical Hamiltonian flows.

Let $(M,g)$ be a Riemannian manifold equipped with a closed two-form
$\sigma$ and a function $V \co M \to \R$ which has a minimum value
of zero. To this data we associate the Hamiltonian flow on the
cotangent bundle $\pi \co T^*M \to M$ defined by the {\em twisted}
symplectic form
$$
\Omega_{\sigma}= d \lambda + \pi^* \sigma
$$
and the {\em total energy} Hamiltonian $H_V \co T^*M \to \R$ given,
in canonical $q,p$--coordinates, by
$$
H_V(q,p) =  \12 g^{-1}(q)(p,p) + V(q).
$$
The flow  of $H_V$ with respect to $\Omega_{\sigma}$ describes  the
classical motion of a charged particle on $M$ under the influence of
the magnetic field $\sigma$ and the (electro-static) force with
 potential $V$. For simplicity we will refer to this as the electro-magnetic
flow for the pair $(\sigma,V)$.

The problem of establishing the existence of periodic orbits for
these flows has been studied extensively. When $V=0$ and $\sigma=0$
this corresponds to the famous problem of finding closed geodesics
for $(M,g)$. When $\sigma=0$, Bolotin proved the existence of closed
orbits on every regular level set in \cite{bo}. When  $V=0$, it is
known from \cite{gi:survey} that there may exist regular energy
levels without periodic orbits. However, for these purely magnetic
flows, it is believed that closed orbits exist on {\em all} small
energy levels. Indeed, this has been established in many cases,
(see, for example, \cite{ar,cmp,gi:surface,ke}).

Schlenk's recent work leads to the following almost existence result
in this setting.
\begin{Theorem}[Schlenk \cite{schl}]
For the electro-magnetic flow determined by a pair $(\sigma,V)$
with $\sigma \neq 0$, there are contractible periodic orbits on
almost every sufficiently small energy level.
\end{Theorem}

From Theorem \ref{thm:ae-area} and Theorem \ref{thm:main} we
immediately obtain the following refinement of this result under
several extra hypotheses.

\begin{Theorem}
\label{thm:em-bound} Consider the electro-magnetic flow determined
by a pair $(\sigma,V)$ where the magnetic two-form is
nondegenerate. Suppose that $(M,\sigma)$ is symplectically
aspherical  and the unit normal bundle of $M$ in $T^*M$ is
homologically trivial in degree $\dim(M)$. Then there is a number
$K>0$ such that almost every sufficiently small energy level
carries a contractible periodic orbit with symplectic area less
than $K$.
\end{Theorem}

\begin{Remark} The conditions on $(M, \sigma)$ are automatically
satisfied in the important case when $(M, \sigma)$ is an
even-dimensional torus with any symplectic form.
\end{Remark}

\subsection{Example: nondegenerate magnetic fields on the flat torus}
\label{flat} According to  Lemma \ref{lem:apriori}, the area bounds
in Theorem \ref{thm:em-bound} can sometimes be used to pass from
almost existence to genuine existence results. To the knowledge of
the author, the only cases where this strategy has been used to
obtain existence results is for classical flows with no magnetic
term ($\sigma=0$), see \cite{hz:wc,hv:cot,hz:book}. Below we show
that this strategy can also be used to partially recover the
existence results of Arnold from \cite{ar} for nondegenerate
magnetic flows on the flat torus.

Let $\T^2$ be the two-dimensional torus
equipped with its standard flat metric. For angular
coordinates  $q_1$ and $q_2$ on the torus, we let $p_1$ and $p_2$
denote the conjugate momenta so that $(q_1,q_2, p_1, p_2)$ are
global coordinates on $T^*\T^2$. We then consider magnetic two-forms
of the form
$$
\sigma =F(q_1,q_2)\,dq_1 \wedge dq_2
$$
where the function $F$ is nonvanishing and positive. The
corresponding Hamiltonian is the kinetic energy
$$
H=  \12 (p_1^2 +p_2^2).
$$
If we set $p_1=r\cos \theta $ and $p_2 = r \sin \theta$, then the
dynamics on the level $\{H=E\}$ is described  by the following
deceptively simple system of equations
\begin{equation}\label{mag}
    \dot{q_1}= \sqrt{2E} \cos \theta,\,\,\,\,\,
  \dot{q_2} = \sqrt{2E} \sin \theta,\,\,\,\,\,
  \dot{\theta}= -F(q_1,q_2).
\end{equation}

The fact that $F$ is nonvanishing can be used to establish the
existence of periodic orbits on fixed energy levels. One way to do
this is to note that if $F$ is nonvanishing, then $\theta$ is
strictly decreasing under the flow. Following \cite{ar}, this allows
one to define, for every energy value $E$,  a Poincar\'{e} return
map $\psi_E \co \T^2 \to \T^2$ whose fixed points correspond to
periodic orbits on the level $\{H=E\}$. In particular, if $\phi^t$
denotes the flow on $\{H=E\}$, then $\psi_E(q_1,q_2)$ is defined by
the equation
$$
\phi^t[(q_1,q_2),0] = [\psi_E(q_1,q_2),-2\pi].
$$
In \cite{ar}, Arnold observes that the return maps $\psi_E$ are all
Hamiltonian diffeomorphisms of $\T^2$, (see \cite{le} for more
details). Applying Conley and Zehnder's proof of the Arnold
Conjecture for symplectic tori \cite{cz}, he then obtains:
\begin{Theorem}[Arnold \cite{ar}]
\label{thm:old-torus} For a nondegenerate magnetic field on the flat
torus there are at least three distinct contractible periodic orbits
on every level set and at least four if they are nondegenerate.
\end{Theorem}

For  nondegenerate magnetic flows in higher dimensions, these return
maps do not exist and there is not yet a way to generalize Arnold's
argument. We now show that it is possible to partially recover
Arnold's result by using the fact that $F$ does not vanish to
establish the {\em a priori} bounds required in Lemma
\ref{lem:apriori}. It is hoped that this strategy can also be used
in higher dimensions.

To begin, we note that when Theorems \ref{thm:main} and
\ref{thm:ae-area} are applied to the magnetic flows above we get:
\begin{Proposition}\label{prop:flat}
For a nondegenerate magnetic two-form $\sigma =F(q_1,q_2)\,dq_1
\wedge dq_2$ on the flat torus there is a $K>0$ such that almost
every low energy level carries a contractible periodic orbit $x$
with symplectic area $|\Aa(x)|<K$.
\end{Proposition}

At this point we observe that since $F$ is nonvanishing, the
projection to $\T^2$ of a closed trajectory of \eqref{mag} is a
closed curve with strictly negative curvature. This observation can
be used to obtain the desired {\em a priori} area bounds when the relative
variance of $F$ is small. More precisely, let $\overline{F}$ and
$\underline{F}$ denote the maximum and minimum values of $F$,
respectively, and set $V_F= \overline{F} / \underline{F}$.
\begin{Proposition} \label{prop:areabound}
If $V_F < \sqrt{\frac{\pi}{2}}$, then there is a continuous function $C(E)>0$
such that
$$
C(E) T \leq |\Aa(x)|.
$$
for every periodic orbit $x$ of \eqref{mag} with period $T$ and
energy $E$.
\end{Proposition}

By Lemma \ref{lem:apriori}, we will then have
\begin{Theorem}
\label{thm:new-torus} Let $\sigma =F(q_1,q_2)\,dq_1 \wedge dq_2$ be
a nondegenerate magnetic two-form on the flat torus such that
$$
V_F <\sqrt{\frac{\pi}{2}}.
$$
For the corresponding magnetic flow, there is a contractible
periodic orbit on every sufficiently low energy level. Moreover, the
symplectic areas of these orbits are uniformly bounded.
\end{Theorem}

\begin{Remark}
The orbits detected in Theorems \ref{thm:old-torus} and
\ref{thm:new-torus} may be different. In particular,  Theorem
\ref{thm:old-torus} finds closed orbits on $\{H=E\}$ which wrap once
around the fibre, whereas the orbits detected by Theorem
\ref{thm:new-torus} are distinguished only by their symplectic area.
\end{Remark}

\begin{proof}[Proof of Proposition~\ref{prop:areabound}]
Let $x(t)$ be a periodic orbit of \eqref{mag} with period $T$ and energy
$E$.  Let $\gamma(t)$ denote the projection of $x(t)$ to $\T^2$.

The symplectic area $\Aa(x)$  can be divided into two terms
\begin{eqnarray*}
\Aa(x) &=& \int_{D^2} \overline{x}^* \Omega_{\sigma} \\
       &=& \int_{D^2} \overline{x}^* d\lambda +  \int_{D^2} \overline{x}^*(\pi^* \sigma)\\
       &=& \int_{D^2} \overline{x}^* d\lambda +  \int_{D^2} \overline{\gamma}^* \sigma\\
       &=& \Aa_1(x) + \Aa_2(\gamma).
 \end{eqnarray*}
A simple computation yields
$$
 \Aa_1(x) = 2E T.
$$
Now let $k=\frac{1}{2\pi} \int_0^T \dot{\theta}(t)\,dt$ be the total
rotation number of $\gamma$. It is clear from \eqref{mag}, that
\begin{equation}\label{area1}
\underline{F}T \leq -2\pi k \leq \overline{F}T.
\end{equation}
Hence,
\begin{equation}\label{area2}
 \Aa_1(x) \geq -k(4E\pi/\overline{F}).
\end{equation}
To bound $\Aa_2(\gamma)$ from below, we note that 
\eqref{mag} implies that the projection  $\gamma(t)$ has constant
speed $\sqrt{2E}$ and negative curvature equal to $-F(\gamma(t))/\sqrt{2E}.$

\begin{Proposition}\label{prop:area}
Let $\xi \co [0,T] \to \R^2$ be a closed planar curve with constant
speed $v$, rotation number $k$, and positive curvature $K(t)/v$ such
that $0 < \underline{K} \leq K(t)$. The Euclidean area enclosed by
$\xi$, $A(\xi)$, satisfies
$$
0 \leq A(\xi) \leq k 4(v/\underline{K})^2.
$$
\end{Proposition}

The proof of this result is elementary and is included in Appendix
A.

Applying Proposition \ref{prop:area} to the curve $\xi(t)
=\gamma(T-t)$ we get the following bounds for the Euclidean area
enclosed by $\gamma$,
$$
0 \geq A(\gamma) \geq k(8E/\underline{F}^2).
$$
For $\Aa_2(\gamma)$, the area of $\gamma$ with respect to $\sigma$,
we then get
\begin{equation}\label{area3}
\Aa_2(\gamma) \geq k\overline{F}(8E/\underline{F}^2).
\end{equation}
Taken together, inequalities \eqref{area1}, \eqref{area2} and
\eqref{area3} imply that
\begin{eqnarray*}
\Aa(x) &\geq& -k\left( \frac{4E\pi}{\overline{F}}
-\frac{8E \overline{F}}{\underline{F}^2} \right)\\
       &=&  -k\frac{8E}{\overline{F}}\left(\frac{\pi}{2} - V_F^2\right)\\
       &\geq& T \frac{4E}{\pi V_F}\left( \frac{\pi}{2} - V_F^2 \right)\\
       &=& T C(E).
\end{eqnarray*}
Clearly, $C(E)$ is positive for $1 \leq V_F < \sqrt{\frac{\pi}{2}}$
and Proposition \ref{prop:areabound} follows.
\end{proof}

\subsection[Hofer's Geometry and c_{HZ}^{1, infty}]{Hofer's Geometry and $c_{\HZ}^{1, \infty}$}
\label{sec:lower}

Theorem \ref{thm:main} can also be used to obtain new lower bounds
for other symplectic invariants of symplectic tubular neighborhoods.
These invariants are defined in terms of Hofer's geometry on
 the group of Hamiltonian diffeomorphisms.

Recall that each smooth, time-dependent, compactly supported
Hamiltonian $H \in C^{\infty}_c([0,1] \times W)$ determines a
Hamiltonian flow, $\phi_H^{t\in [0,1]}$, on $(W, \Om)$. The group of
Hamiltonian diffeomorphisms of $(W, \Om)$, $\Ham(W, \Om)$, consists
of all the time--1 maps $\phi_H^1$ obtained in this manner.
Conversely, for every path $\phi^{t \in [0,1]}$ in $\Ham(W, \Om)$
there is a Hamiltonian $H \in C^{\infty}_c([0,1] \times W)$ which
generates it, ie $\phi^t = \phi_H^t \circ \phi^0.$ When $(W^{2l},
\Om)$ is closed, this generating Hamiltonian is uniquely determined
by the normalization condition
$$
\int _W H(t,x)\, \Om^{l} = 0, \text{ for all  }t \in [0,1].
$$
Otherwise, the generating Hamiltonian is already uniquely determined
by the fact that it has compact support.

In \cite{ho1}, Hofer defines the length of a Hamiltonian path
$\phi^t$, in terms of its generating Hamiltonian $H$, by the formula
\begin{eqnarray*}
\Length(\phi^t) &=& \int_0^1 \max_{x \in W}H(t,x)\, dt - \int_0^1\min_{x \in W}H(t,x)\, dt\\
               &=& \Length^+(\phi^t) +  \Length^-(\phi^t).
\end{eqnarray*}
The Hofer metric on $\Ham(W, \Om)$ is then defined by
$$
\rho(\phi, \psi) = \inf \Length(\phi^t),
$$
where the infimum is taken over all smooth paths joining $\phi$ to
$\psi$.  A path $\phi^{t \in [0,1]}$ is said to be \emph{length
minimizing} if $\Length(\phi^t)=\rho(\id,\phi^1)$.

The quantity $\Length^+(\phi^t)$ is called the positive Hofer length
of $\phi^t$, and one defines the positive Hofer metric $\rho^+$ in
the obvious way.

One of the central problems in the study of Hofer's geometry is to
characterize the Hamiltonian paths which minimize the (positive)
Hofer length. This question is profoundly related to the dynamical
features of the corresponding flow \cite{pol, po:book}. It is of
particular interest for paths which are generated by
time-independent Hamiltonians (referred to here as autonomous
Hamiltonian paths), because these describe important dynamical
systems from classical mechanics. The first result in this direction
was proved by Hofer in \cite{ho2}. There he shows that an autonomous
Hamiltonian path in $\Ham(\R^{2l},\Om_{2l})$ is length minimizing as
long as its flow has no nonconstant periodic orbits with period less
than or equal to one. This result has been extended in several
directions \cite{bp,kl,lm,mcsl, mc:geom,oh,sib,us} and Polterovich
has conjectured that it is true for general symplectic manifolds
\cite{po:book}.

In \cite{sik}, Sikorav shows that the length minimizing property of
an autonomous Hamiltonian path is also constrained by the set on
which it is supported. This motivates the definition of the
following notion, which we call the \emph{Hofer width} of an open
subset $U$ of $(W, \Omega)$:
$$
w_H(U) \eqdef \sup\left\{ \max(H)-\min(H) \mid H \in
C^{\infty}_c(U),\,\,\Length(\phi_H^t)=\rho(\id ,\phi_H^1)\right\}.
$$
This invariant measures the size of $U$ in $W$ in terms of the
length of the longest autonomous Hamiltonian path which is supported
in $U$ and is length minimizing in $\Ham(W, \Om)$.

The positive Hofer width  $w^+_H(U)$, is defined analogously as
follows
$$
w_H^+(U) \eqdef \sup\left\{ \max(H) \mid H \in
C^{\infty}_c(U),\,H\geq 0,\,\,\Length^+(\phi_H^t)=\rho^+(\id
,\phi_H^1)\right\}.
$$
When $W$ is noncompact, the positive Hofer width $w_H^+(U)$ can be
interpreted in the same way as $w_H(U)$. When $W$ is closed this
interpretation must be altered slightly, because the quantity
$\max(H)$ is no longer equal to the positive length of $\phi_H^t$.
In particular, each  $H \in C^{\infty}_c(W)$ is not properly
normalized when $W$ is closed, and
\begin{equation} \label{positive}
\Length^+(\phi_H^t) = \max(H) - \frac{1}{\volume(W)}\int_W H(t,x)\,
\Om^{l}.
\end{equation}
Regardless of this, it is evident from \eqref{positive} that
$\max(H)$ is a reasonable substitute for the positive Hofer length
since
$$
\max(H) \geq \Length^+(\phi_H^t) \geq \max(H)\left( 1 - \volume(U)/
\volume(W) \right).
$$
For a closed symplectic manifold  $(W, \Om)$, the work of McDuff 
\cite{mc:geom} implies that the (positive) Hofer width of any open
subset $U \subset W$ is nonzero. For tubular neighborhoods of
certain submanifolds much more can be said. For example, suppose
that in addition to being closed, $(W, \Om)$ is symplectically
aspherical. Let $L \subset W$ be a Lagrangian submanifold such that
the map $\pi_1(L) \to \pi_1(W)$ (induced by inclusion) is injective,
and $L$ admits a metric $g$ with no contractible geodesics, eg
$L$ is a torus or admits a hyperbolic metric. In this case, it
follows from the work of Schwarz \cite{sc} that the (positive)
Hofer width of any neighborhood of $L$ is infinite. In particular,
it is easy to construct a positive autonomous Hamiltonian $H$
supported in any neighborhood of $L$, such that $\max(H)$ is
arbitrarily large and the flow of $H$ has no nonconstant
contractible periodic orbits of any period. Schwarz proves in
\cite{sc} that the latter property implies that $\phi_H^t$ minimizes
the (positive) Hofer length.

In \cite{kl}, this sufficient condition for a Hamiltonian path
to be length minimizing is refined. The following result can be easily
derived from \cite[Theorem~1.8]{kl}.

\begin{Proposition}[Kerman--Lalonde \cite{kl}]
Let $H$ belong to $\Hh_{\test}(U)$ and suppose that the flow of $H$
has no nonconstant contractible periodic orbits with period equal to
one and action greater than $\max(H)$. Then the path
$\{\phi_H^t\}_{t \in [0,1]}$ minimizes the positive Hofer length.
\end{Proposition}

In other words, every admissible test Hamiltonian used to define
$c_{\HZ}^{1,\infty}(U)$ generates a Hamiltonian path which minimizes
the positive Hofer length. From this, we immediately get the
following inequality.
\begin{Proposition}
\label{prop:ineq} Suppose that  $(W, \Om)$ is closed and
symplectically aspherical. For every open set $U$ in $W$ we have
$$
w^+_H(U) \geq c_{\HZ}^{1,\infty}(U).
$$
\end{Proposition}

Theorem \ref{thm:main} now implies the following result.

\begin{Theorem}
For any symplectic tubular neighborhood $U_R$ as in Theorem
\ref{thm:main}, we have $$w^+_H(U_R) \geq \pi R^2.$$
\end{Theorem}

\begin{Remark}
Instead of the Hofer width of $U$, one may consider the extrinsic
diameter of $\Ham (U, \Om)$, ie the diameter of $\Ham (U, \Om)$
as a subset of $\Ham (W, \Om)$. This is clearly bounded from below by the
positive Hofer width. However, Ostrover \cite{os} has recently  proven that if
$W$ is closed and $\pi_2(W)=0$, then the extrinsic diameter of
$\Ham(U,\Om)$ is infinite for any open set $U \subset W$.\footnote{The
author thanks the referee for informing him of this point.}
The intrinsic diameter of $\Ham (U, \Om)$ is unlikely to be
interesting, as Polterovich has conjectured that the diameter of
$\Ham(W, \Om)$ for any symplectic manifold is infinite \cite{pol}.
\end{Remark}

\begin{Remark}
Using Proposition \ref{prop:ineq}, we can generalize the discussion
above concerning the Hofer width of neighborhoods of Lagrangian
submanifolds. Let $L$ be a Lagrangian submanifold of $(W, \Om)$
which we again assume is closed and symplectically aspherical. Let
$g$ be any metric on $L$. By Weinstein's Neighborhood theorem there
is a neighborhood $V_R$ of $L$ in $W$ which is symplectomorphic to a
neighborhood of the zero-section in $(T^* L, d\lambda)$ of the form
$\{ (q,p) \in T^* L \mid \|p\|_g^2=R^2\}$. Set $T(g,L,W)$ equal to
the minimal period of all closed geodesics of $g$ which represent a
homotopy class in the kernel of the map $\pi_1(L) \to \pi_1(W).$ A
simple argument shows that $ c_{\HZ}^{1,\,\infty}(V_R) \geq T(g,L,W)
\, R.$ Hence, by Proposition \ref{prop:ineq}, $w_H^+(V_R) \geq
T(g,L,W) \, R$.
\end{Remark}

\section{Symplectic tubular neighborhoods}
\label{sec:U_R}

In this section we define the symplectic tubular neighborhoods $U_R$
of $M$ in $(W, \Om)$. We also show that there are Hamiltonian flows
on each $U_R$ which are totally periodic. These flows are then used
to obtain lower bounds for the capacity of $U_R$. This allows us to
reduce the proof of Theorem \ref{thm:main} to a statement concerning
the dynamics of test Hamiltonians on symplectic tubular
neighborhoods.

\subsection{The definition}

Let $\pi \co E\to M$ be the symplectic normal bundle to $M$, ie
the symplectic orthogonal complement to $TM$ in $TW|_M$. Henceforth,
we will identify $M$ with the zero-section of $E$ and write points
in $E$ as pairs $(p,z)$, where $p \in M$ and $z$ is in the fibre of
$E$ over $p$. The bundle $E$ is a symplectic vector bundle with a
fibrewise symplectic form $\sigma_E$ given by the restriction of
$\Omega$ to the fibres. Let $J_E$ be a fibrewise almost complex
structure on $E$ which is compatible with $\sigma_E$. This yields a
fibrewise inner product $g_E (\cdot,\cdot)=\sigma_E(\cdot, J_E
\cdot)$ and we denote the function $(p,z) \mapsto g_E(p)(z,z)$ by
the (over-simplified) notation $\|z\|^2$. For these choices, we set
$$
U_R=\{ (p,z) \in E \mid \|z\|^2  \leq R^2 \}.
$$
On $E \setminus M$ there is a canonical fibrewise one-form
$\alpha_E$ which is defined by the formula
$$
\alpha_E(p,z)(\cdot) = \|z\|^{-2} \sigma_E(p)(z,\cdot)
$$
and satisfies
$$
d^F(\12 \|z\|^2 \alpha_E) = \sigma_E,
$$
where $d^F$ denotes the fibrewise exterior differential. Fixing a
Hermitian connection $\nabla$ on $E$ we can extend $\alpha_E$ to a
genuine one-form on $E \setminus M$ which we denote by $\alpha$. We
now consider the closed two-form
$$
d(\12 \|z\|^2 \alpha)+\omega
$$
which is defined on all of $E$ and is nondegenerate, and hence
symplectic, on a neighborhood of the zero section. Here we have
identified $\omega=\Omega|_M$ with its pull-back to the total space
$E$.

By Weinstein's Symplectic Neighborhood Theorem, $(U_R, d(\12 \|z\|^2
\alpha) + \om)$ is symplectomorphic to a neighborhood of $M$ in
$(W,\Om)$ when $R>0$ is sufficiently small. We will restrict our
attention to values of $R$ for which this symplectomorphism
exists.\footnote{ This is the first condition restricting the size
of $R$.} We refer to these neighborhoods as \emph{symplectic tubular
neighborhoods} of $M$.

\subsection{Dynamics on $U_R$}

\begin{Lemma}
\label{periodic} On each $U_R$, the Hamiltonian flow of the
function $\|z\|^2$ with respect to the symplectic form $\Omega= d(\12 \|z\|^2
\alpha) + \om$ is periodic with period $\pi$.
\end{Lemma}

\begin{proof}
Let $V_{\Hopf}(p,z)= J_E(p)z$ be the vector field on $E$ which generates the
standard fibrewise Hopf circle action with period $2\pi$.  Then
\begin{eqnarray*}
  i_{V_{\Hopf}}(d (\12 \|z\|^2 \alpha) + \om)&=& i_{V_{\Hopf}}(d (\12 \|z\|^2 \alpha)) \\
    &=& \Ll_{V_{\Hopf}} (\12 \|z\|^2 \alpha) - d (i_{V_{\Hopf}}(\12 \|z\|^2 \alpha)) \\
    &=& -d(\12 \|z\|^2),
\end{eqnarray*}
which proves the result.
\end{proof}

\begin{Corollary}\label{cor:lowerbound}
The capacity  $c_{\HZ}(U_R)$ is no less than $\pi R^2$.
\end{Corollary}

\begin{proof}
Let
$h_{\epsilon} \co [0,R^2] \to [0,\, \infty)$ be a smooth function satisfying
\begin{itemize}
    \item $h_{\epsilon}(s) = \pi R^2 - \epsilon$ for $s$ near $0$,
    \item $h_{\epsilon}(s) = 0 $ for $s$ near $R^2$,
    \item $|h_{\epsilon}'|< \pi$ for all $s \in [0,R^2].$
\end{itemize}
Consider the test Hamiltonian $H_{\epsilon}(p,z) =
h_{\epsilon}(\|z\|^2)$ on $U_R$.
The Hamiltonian vector field of $H_{\epsilon}$ is $X_{H_{\epsilon}}=
2h_{\epsilon}'(\|z\|^2)\cdot V_{\Hopf}$. Since  $|h_{\epsilon}'|<
\pi$, it follows that $H_{\epsilon}$ is admissible for all
$\epsilon>0$.
\end{proof}

The previous corollary also implies that $\smash{c_{\HZ}^{\,1,\,2\pi
R^2}(U_R)} \geq \pi R^2$, and so to prove Theorem \ref{thm:main} we
must show that $\smash{c_{\HZ}^{\,1,\,2 \pi R^2}(U_R)} \leq \pi
R^2$. This upper bound for the refined capacity is implied by the
following result concerning the dynamics of test Hamiltonians on
$U_R$.

\begin{Theorem}\label{thm:main-1dyn}
Suppose that $(W, \Omega)$ is geometrically bounded and
symplectically aspherical. Let $M^{2m}$ be a closed symplectic
submanifold of $W$ whose unit normal bundle is homologically trivial
in degree $2m$. For all sufficiently small $R>0$, if $H \in
\Hh_{\test}(U_R)$ satisfies $\max(H) > \pi R^2$, then $H$ has a
nonconstant $1$--periodic orbit with action in the interval
$(\max(H), \max(H) + 2\pi R^2)$.
\end{Theorem}

The rest of this paper is devoted to the proof of
Theorem \ref{thm:main-1dyn}.

\section{Tools from Floer theory}

To prove Theorem \ref{thm:main-1dyn}, we will use machinery from
Floer theory which takes advantage of the filtration by the action
functional. In this section we recall the definitions and properties
of the necessary tools when $(W,\Om)$ is geometrically bounded and
symplectically aspherical. For such manifolds this material is
essentially standard by now and we recall the required results
without proof. The reader is referred to the sources
\cite{cfh,fl1,fl2,fl3,fh,hs,po,sa} for more details.

In fact, we make use of the Morse complex, the Floer complex for
Hamiltonian diffeomorphisms, and the Floer complex for pairs of
Lagrangians. We begin by describing why all three of these objects
are used in the proof.

Recall that a Morse complex can be associated to a generic choice of
a Morse function $h$ and a metric $g$ on a closed manifold. The
Morse chain complex $(C(h), \p_g)$ consists of the vector space
$C(h)$ generated by the critical points of $h$ and the boundary map
$\p_g$ which counts trajectories of the negative gradient flow of
$h$ with respect to $g$. The homology of this complex is independent
of the choice of the data $h$ and $g$ and is isomorphic to the
homology of the manifold on which they are defined.

Let $\Hh = C_c^{\infty}(S^1 \times W)$ be the space of smooth,
time--periodic, compactly supported functions on a symplectic
manifold $(W, \Om)$. Every $G \in \Hh$ determines an action
functional on the space of smooth contractible loops in $W$. As
well, each $S^1$--family of $\Om$--compatible almost complex
structures, $J_t$, determines a metric on this loop space. The
Hamiltonian Floer complex for a  generic pair $(G, J_t)$ is then
defined, at least heuristically, as the Morse complex for this
functional and metric on the loop space.

Alternatively, the same data, $G$ and $J_t$, determines
 a $[0,\, 1/2]--$family of almost complex structures $\wtilde{J}^W_t$ on  the product manifold
$$
(\wtilde{W},\, \wtilde{\Om}) \eqdef (W \times W,\, \Om \oplus
(-\Om)).
$$
This family can then be used to define a Lagrangian Floer complex
 for the pair of Lagrangian submanifolds of $\widetilde{W}$ consisting of the diagonal
 $\Delta$ and the graph of $\phi^1_G$.

These two Floer complexes associated to $G$ and $J_t$ are identical
(see  Section~\ref{translation}). However, in the Lagrangian setting one
has the added freedom to compute the homology by replacing the graph
of $\phi_G^1$ and $\wtilde{J}^W_t$ by a Lagrangian submanifold and
family of almost complex structures which aren't necessarily
determined by objects on $W$. This freedom will allow us to make
some crucial identifications of certain restricted Lagrangian Floer
complexes with Morse complexes. This is achieved using the
Morse--Bott version of Lagrangian Floer theory developed by Po\'zniak
in \cite{po}, (see Theorem~\ref{thm:poz}). It is not clear whether
these identifications can be obtained using Hamiltonian Floer theory
alone (cf \cite{bps,cfhw}).

In the first part of this section we recall the necessary material
from Lagrangian Floer theory, including Po\'zniak's theory of clean
intersections from \cite{po}. We also discuss the definitions and
properties of monotone Floer continuation maps. In the second part
of the section we recall the basic construction of the Hamiltonian
Floer complex and we describe the natural translation from
Hamiltonian Floer theory to Lagrangian Floer theory. Finally, we
describe some useful algebraic relations satisfied by monotone Floer
continuation maps in the  Hamiltonian setting.

\noindent\textbf{Concerning notation}\qua In what follows, objects
defined on (for) the manifold $\widetilde{W}$ will be decorated with
a tilde. If, in addition, they are defined using an object coming
from $W$, then they will also be given a superscript $W$. Generic
approximations of any object will be denoted by a prime.

\subsection{Lagrangian Floer homology}
\label{sec:fc} We describe here the simplest version of Lagrangian
Floer theory which is defined for a fixed Lagrangian submanifold and
its images under Hamiltonian diffeomorphisms. We will restrict our
attention to the Lagrangian submanifold given by the diagonal
$\Delta$ in the product manifold $(\wtilde{W},\,\wtilde{\Om})= (W
\times W,\, \Om \oplus (-\Om))$.

It is clear from the definitions that if $(W, \Om)$ is
symplectically aspherical and geometrically bounded with respect to
the almost complex structure $J_{gb}$, then
$\smash{\bigl(\wtilde{W},\,\wtilde{\Om}\bigr)}$ is also symplectically aspherical and
is geometrically bounded with respect to the almost complex
structure $\smash{\wtilde{J}_{gb}} = J_{gb} \oplus -J_{gb}$. Let
$$\widetilde{\Hh} = C^{\infty}_c\bigl(S^1 \times \wtilde{W}, \R\bigr)$$
be the space of smooth time--periodic Hamiltonians on $\wtilde{W}$ which
have compact support. Given a $\widetilde{G} \in \widetilde{\Hh}$,
we define  $\smash{\widetilde{\Jj}_{[0, \, 1/2]} =\widetilde{\Jj}_{[0, \,
1/2]}(\wtilde{G})}$ to be the set of all
 smooth $[0, \, 1/2]$--families of $\widetilde{\Om}$--tame almost complex structures
$\widetilde{\Jj}_t$ on $\widetilde{W}$ which are
$\widetilde{\Om}$--compatible near the support of $\wtilde{G}$ and
are equal to $\wtilde{J}_{gb}$ outside some compact set. For a
generic pair $\smash{\bigl(\widetilde{G}, \widetilde{J}_t\bigr) \in \widetilde{\Hh}
 \times \widetilde{\Jj}_{[0, \, 1/2]}}$,  the corresponding
Lagrangian Floer chain complex  for actions in the positive interval
$(a,b)$ is denoted by
 $$
 (\CF^{a,b} (\Delta,\widetilde{G}),
 \partial_{\widetilde{J}_t}),
 $$
and is constructed as follows.

Let $\widetilde{\Ll} = \widetilde{\Ll}(\widetilde{G})$ denote the
space of smooth paths
$$
\big\{ \widetilde{x} \in C^{\infty}([0, \, 1/2],\wtilde{W}) \mid
\widetilde{x}(0) \in \Delta,\, \widetilde{x}(1/2) \in
(\phi_{\widetilde{G}}^1)^{-1}(\Delta) \big\}.
$$
Consider the action functional $\widetilde{\Aa}_{\widetilde{G}} \co
\wtilde{\Ll} \to \R$ defined by
\begin{equation*}\label{lag}
    \widetilde{\Aa}_{\widetilde{G}}(\widetilde{x})= - \int_{[0, \, 1/2]^2}
\widetilde{v}^*
    \wtilde{\Om},
\end{equation*}
where  $\widetilde{v}$ is a map from  $[0, \, 1/2] \times [0, \, 1/2]$
to $\wtilde{W}$ such that $\widetilde{v}(s,\cdot) \in \wtilde{\Ll}$ for
all $s \in [0, \, 1/2]$, $\widetilde{v}(1/2,t)= \widetilde{x}(t)$,  and
$\widetilde{v}(0,t) = \gamma_0(t)$ for some fixed $\gamma_0 \in
\wtilde{\Ll}.$
The critical point set of $\widetilde{\Aa}_{\widetilde{G}}$,
$\Crit(\widetilde{\Aa}_{\widetilde{G}})$,  consists of the constant
paths in $\widetilde{\Ll}$ and hence coincides with the intersection
points of the diagonal $\Delta$ and its image
$(\phi_{\widetilde{G}}^1)^{-1}(\Delta)$. A critical point
$\widetilde{x}$ is nondegenerate if the corresponding point of
intersection is transverse. In this case, $\widetilde{x}$ has a
well-defined Maslov index, $\mu_{\Maslov}(\widetilde{x},\widetilde{G})$.

For constants $0<a<b$, let $\widetilde{\Hh}^{a,b} \subset
\widetilde{\Hh}$ be the open set of all functions $\widetilde{G}$ in
$\widetilde{\Hh}$ such that  $a$ and $b$ are not critical values of
$\widetilde{\Aa}_{\widetilde{G}}$. For $\widetilde{G} \in
\widetilde{\Hh}^{a,b}$, let
$\smash{\Crit^{a,b}(\widetilde{\Aa}_{\widetilde{G}})}$ be the set of
critical points of $\smash{\widetilde{\Aa}_{\widetilde{G}}}$ with action in
the interval $(a,b)$. Generically, the elements of
$\smash{\Crit^{a,b}(\widetilde{\Aa}_{\widetilde{G}})}$ are nondegenerate and
hence finite in number. In this case, the Lagrangian Floer chain
group of  $\widetilde{G}$ for actions restricted to $(a,b)$ is
defined to be the graded $\Z_2$--vector space
$$
\CF^{a,b} (\Delta, \widetilde{G})= \bigoplus_{\widetilde{x} \in
\Crit^{a,b}(\widetilde{\Aa}_{\widetilde{G}})} \Z_2 \widetilde{x}.
$$
The Floer boundary operator $\p_{\widetilde{J}_t}$ is defined using
the  Floer moduli spaces
$\smash{\wtilde{\Mm}(\widetilde{x},\widetilde{y},\widetilde{G},
\wtilde{J}_t)}$ for
pairs $\smash{\widetilde{x}, \widetilde{y} \in
\Crit^{a,b}(\widetilde{\Aa}_{\widetilde{G}})}$. Each
$\smash{\wtilde{\Mm}(\widetilde{x},\widetilde{y},\widetilde{G},\wtilde{J}_t)}$
consists of the maps $\widetilde{u} \co \R \times [0, \, 1/2] \to
\wtilde{W}$ which are solutions of the equation
\begin{equation}
\label{eq:moduli-cr}
\partial_s \widetilde{u} + \wtilde{J}_{t}(\widetilde{u}) \partial_t\widetilde{u}=0,
\end{equation}
satisfy the boundary conditions
\begin{equation}
\label{eq:moduli-boundary}
\widetilde{u}(s,0) \in \Delta\qquad \text{and}\qquad \widetilde{u}(s,1/2) \in (\phi^1_{\widetilde{G}})^{-1}(\Delta),
\end{equation}
and have limits
\begin{equation}
\label{eq:moduli-limits} \lim_{s \to -\infty}\widetilde{u}(s,t) =
\widetilde{x}(t)\qquad  \text{and}\qquad
 \lim_{s \to \infty}\widetilde{u}(s,t) =
\widetilde{y}(t),
\end{equation}
which are uniform in $t$.
For generic data $\widetilde{G}$ and $\wtilde{J}_t$, every moduli
space $\wtilde{\Mm}(\widetilde{x},\widetilde{y},\widetilde{G}
,\wtilde{J}_t)$ is a smooth manifold of dimension
$\mu_{\Maslov}(\widetilde{x}, \widetilde{G}) - \mu_{\Maslov}(\widetilde{y},
\widetilde{G})$. (If $\widetilde{G}$ is fixed, then any family
$\wtilde{J}_t$ which forms a generic data pair with $\widetilde{G}$
will be said to be \emph{regular for} $\widetilde{G}$.) We also note
that each
$\smash{\wtilde{\Mm}(\widetilde{x},\widetilde{y},\widetilde{G}
,\wtilde{J}_t)}$ is equipped with a free $\R$--action given by $\tau
\cdot \widetilde{u}(s,t)= \widetilde{u}(s + \tau, t)$.

Let $\Crit^{a,b}_k(\widetilde{\Aa}_{\widetilde{G}})$ denote the
subset of $\Crit^{a,b}(\widetilde{\Aa}_{\widetilde{G}})$ which
consists of critical points with Maslov index equal to $k$.
The boundary operator
$$\partial_{\widetilde{J}_t} \co
  \CF^{a,b}(\Delta, \widetilde{G}) \to \CF^{a,b}(\Delta,
  \widetilde{G})$$
is then defined on each $\widetilde{x} \in
\smash{\Crit^{a,b}(\widetilde{\Aa}_{\widetilde{G}})}$ by the formula
$$
\partial_{\widetilde{J}_t} (\widetilde{x}) = \sum_{\scriptsize
  \begin{array}{c}
  \widetilde{y} \in \Crit^{a,b}(\widetilde{\Aa}_{\widetilde{G}}),\\
  \mu_{\Maslov}(\widetilde{y},\widetilde{G})
  = \mu_{\Maslov}(\widetilde{y},\widetilde{G})-1
  \end{array}} \hspace{-40pt}\#
(\widetilde{\Mm}(\widetilde{x},\widetilde{y},\widetilde{G},\widetilde{J}_t)
/ \R) \widetilde{y},
$$
where $\#
(\widetilde{\Mm}(\widetilde{x},\widetilde{y},\widetilde{G},\widetilde{J}_t)
/ \R)$ is the number  of elements in the $0$--dimensional manifold
$\Mm(\widetilde{x},\widetilde{y}, \widetilde{G}, \widetilde{J}_t)/\R$,
modulo $2$. Since we are assuming that $(W, \Om)$ is symplectically
aspherical, no bubbling can occur and there are no obstructions. So,
the usual arguments imply that
 $$
\partial_{\widetilde{J}_t} \circ \partial_{\widetilde{J}_t} = 0.
$$
The corresponding homology groups $\HF^{a,b}(\Delta,
\widetilde{G})$ are independent of the choice of the regular family
$\widetilde{J}_t$. They are also locally constant on
$\widetilde{\Hh}^{a,b}$. This allows one to define $\HF^{a,b}
(\Delta, \widetilde{G})$ for any $\widetilde{G} \in
\widetilde{\Hh}^{a,b}$, regardless of whether or not the elements of
$\Crit^{a,b}(\widetilde{\Aa}_{\widetilde{G}})$ are nondegenerate.
One just sets $\HF^{a,b} (\Delta, \widetilde{G})=\HF^{a,b}(\Delta,
\widetilde{G}')$ for some $\widetilde{G}'  \in
\widetilde{\Hh}^{a,b}$ which is close to $\widetilde{G}$ and has the
required nondegeneracy property.

\subsection{Clean intersections}
\label{sec:ci}

The extension of the definition of the restricted Floer homology
$\HF^{a,b} (\Delta, \widetilde{G})$ to every $\widetilde{G} \in
\widetilde{\Hh}^{a,b}$ is particularly useful in the case when
$\Delta$ and $(\phi^1_{\widetilde{G}})^{-1}(\Delta)$ intersect
nicely along submanifolds. This situation was first studied by
Po\'zniak in \cite{po} (for more general versions of Lagrangian Floer
homology). In this section, we briefly describe Po\'zniak's results
and  some refinements from \cite{bps}, as they apply to our present
setting.

A submanifold $\widetilde{N} \subset \Delta \cap
(\smash{\phi_{\widetilde{G}}^1})^{-1}(\Delta)$ is said to be a \emph{clean
intersection} for $\widetilde{G} \in \widetilde{\Hh}$ if it is a
connected component of $\Delta \cap
(\phi^1_{\widetilde{G}})^{-1}(\Delta)$ and
$$
T_p (\widetilde{N}) = T_p (\Delta) \cap T_p(
(\phi^1_{\widetilde{G}})^{-1}(\Delta)) \text{    for all   } p \in
\widetilde{N}.
$$
For clean intersections, Po\'zniak proved the following result.

\begin{Theorem}{\rm\cite[Theorem 3.4.11]{po}, \cite[Theorem
5.2.2]{bps}}\qua
\label{thm:poz} Suppose that the set $\widetilde{N}$ is a clean
intersection for $\widetilde{G} \in \widetilde{\mathcal{H}}^{a,b}$
and that $\widetilde{N} = \Crit^{a,b}
(\widetilde{\Aa}_{\widetilde{G}})$. Let $h \co \widetilde{N} \to \R$
be a Morse function and let $g$ be a metric on $\widetilde{N}$ such
that the corresponding Morse complex $(C(h),\partial_g)$ is
well-defined. Then there is a Hamiltonian $\widetilde{H}_h \in
\widetilde{\Hh}$, a family $\widetilde{J}_{g,t} \in
\widetilde{\Jj}_{[0,\,1/2]}$, and a constant $\delta_0
>0$ such that for every $\delta < \delta_0$ we have the following strong
equivalence of complexes
\begin{equation}\label{eq:complex-id}
(\CF^{a,b} (\Delta, \widetilde{G}+ \delta \widetilde{H}_h) ,
 \partial_{\widetilde{J}_{g,t}}) \equiv
 (C(h),\partial_g).
\end{equation}
\end{Theorem}

The equivalence relation ``$ \equiv $" denotes the fact that the
elements of the set  $\smash{\Crit^{a,b}(\widetilde{\Aa}_{\widetilde{G}+
\delta \widetilde{H}_h})}$ actually coincide with the critical points
of $h$ and the elements of the one--dimensional moduli spaces
$\wtilde{\Mm}(\widetilde{x},\widetilde{y},\widetilde{G}+ \delta
\widetilde{H}_h,\wtilde{J}_{g,t})$ are $t$--independent and coincide
with the negative gradient trajectories of $h$ with respect to $g$.

It follows from \eqref{eq:complex-id} that
\begin{equation}\label{eq:hfl=m}
\HF^{a,b}_* (\Delta, \widetilde{G}) =
H_{*-\mu_{\Maslov}(\widetilde{N}, \widetilde{G})}
(\widetilde{N};\Z_2).
\end{equation}
We will refer to the grading shift $\mu_{\Maslov}(\widetilde{N} ,
\widetilde{G})$ as the \emph{relative Maslov index} of
$\widetilde{N}$. It depends on the linearized flow of
$\widetilde{G}$ along $\widetilde{N}$.

\subsubsection{Outline of the proof of Theorem~\ref{thm:poz}}

To establish the strong equivalence of complexes in equation
\eqref{eq:complex-id}, Po\'zniak introduces the notion of a local
Floer complex of $\widetilde{N}$ and proves that there is a local
complex which is strongly equivalent to both sides.

To describe this argument, we first require the notion of a
$\widetilde{J}_t$--isolating neighborhood of $\widetilde{N}$ in
$\widetilde{W}$. This is a precompact neighborhood $\widetilde{U}$
of $\widetilde{N}$ such that any solution $\widetilde{u} \co \R \times
[0, 1/2 ] \to \widetilde{W}$ of \eqref{eq:moduli-cr} and
\eqref{eq:moduli-boundary} whose image lies in the closure of
$\widetilde{U}$, must satisfy $\widetilde{u}(s,t) = \widetilde{x} \in
\widetilde{N}$ for all $(s,t) \in \R \times [0, 1/2 ]$. Every clean
intersection admits a $\widetilde{J}_t$--isolating neighborhood for
any choice of the family $\widetilde{J}_t$ (see \cite[Lemma
5.2.3]{bps}).

Let $\widetilde{U}$ be a $\widetilde{J}_t$--isolating neighborhood
for $\widetilde{N}$ and choose a  $\widetilde{G}' \in
\widetilde{\Hh}$ such that $\widetilde{G}'$ is $C^2$--close to
$\widetilde{G}$ and the
 elements of $\Crit(\widetilde{\Aa}_{\widetilde{G}'})$ are nondegenerate. The local
Floer chain complex of $\widetilde{G}'$ is defined to be the
$\Z_2$--vector space
\begin{equation*}
    \CF^{\loc}(\widetilde{G}',\widetilde{U}) = \bigoplus_{\widetilde{x}' \in
    \Crit(\widetilde{\Aa}_{\widetilde{G}'}), \, \widetilde{x}' \subset
    \widetilde{U}}\Z_2 \widetilde{x}'.
\end{equation*}
For a family $\widetilde{J}_t'$ which is regular for
$\widetilde{G}'$ and is $C^1$--close to $\widetilde{J}_t$, the
boundary operator $\partial^{\loc}_{\widetilde{J}_t'}$ for the
complex is defined in the usual way except that it only counts
solutions $\widetilde{u}$ of \eqref{eq:moduli-cr} whose images lies in
$\widetilde{U}$. The corresponding homology is called the local
Floer homology of $\widetilde{N}$,
$$
\HF^{\loc}(\widetilde{N}) \eqdef H_*(
\CF^{\loc}(\widetilde{G}',\widetilde{U}),
\partial^{\loc}_{\widetilde{J}_t'}),
$$
and Po\'zniak proves that it is independent of the choices of
$\widetilde{U}$, $\widetilde{G}'$ and $\widetilde{J}_t'$.

To calculate $\HF^{\loc}(\widetilde{N})$, Po\'zniak starts with a
Morse function $h \co \widetilde{N} \to \R$ and a generic metric $g$
on $\widetilde{N}$ and he constructs explicit approximations of
$\widetilde{G}$ and $\widetilde{J}_t$ whose local Floer complex is
strongly equivalent to the Morse complex determined
 by $h$ and $g$. More precisely, Po\'zniak proves
that there is a Hamiltonian $\widetilde{H}_h \in \widetilde{\Hh}$, a
family of almost complex structures $\widetilde{J}_{g,t}$, and a
$\delta_0>0$ such that for $\delta < \delta_0$ one has
$$
(\CF^{\loc}(\widetilde{G}+\delta
\widetilde{H}_h,\widetilde{U}),
\partial^{\loc}_{\widetilde{J}_{g,t}}) \equiv (C(h),\partial_g)
$$
for any choice of $\widetilde{U}$ (see \cite[Proposition
3.4.6]{po}).

To complete the proof of Theorem \ref{thm:poz}, it only remains to show
that if
$\widetilde{N}=\Crit^{a,b}(\widetilde{\Aa}_{\widetilde{G}})$ then
$$
(\CF^{\loc}(\widetilde{G}+\delta
\widetilde{H}_h,\widetilde{U}),
 \p^{loc}_{\widetilde{J}_{g,t}} )
\equiv ( \CF^{a,b}(\Delta, \widetilde{G} + \delta
\widetilde{H}_h), \p_{\widetilde{J}_{g,t}} )
$$
for sufficiently small $\delta>0$. This follows immediately from
\cite[Lemma 5.2.4]{bps}.

\subsection{Monotone Floer continuation maps}
\label{sec:mon-hom}

As mentioned above, the restricted Floer homology $\HF^{a,b}(\Delta,
\widetilde{G})$ is only locally constant on $\widetilde{\Hh}^{a,b}$.
However, if the functions $\widetilde{G},\widetilde{H} \in
\widetilde{\Hh}^{a,b}$ satisfy $\widetilde{G} \geq \widetilde{H}$,
ie $\widetilde{G}(t,x) \geq \widetilde{H}(t,x)$ for all $(t,x)
\in S^1 \times \widetilde{W}$, then there is a natural homomorphism
$$
\sigma_{\widetilde{H}\widetilde{G}} \co \HF^{a,b} (\Delta,
\widetilde{G}) \to \HF^{a,b}(\Delta,\widetilde{H})
$$
called a \emph{monotone Floer continuation map}. Since our later
arguments rely heavily on these maps, we recall their definition.

For $\widetilde{G},\widetilde{H} \in \widetilde{\mathcal{H}}^{a,b}$
with $\widetilde{G} \geq \widetilde{H}$, a monotone homotopy from
$\widetilde{G}$ to $\widetilde{H}$ is a family of functions
$\widetilde{G}_s \in \widetilde{\Hh}$ such that $\partial_s
\widetilde{G}_s \leq 0$ and
$$
\widetilde{G}_s=\begin{cases}
\widetilde{G}&\text{for $s\in (-\infty,-1]$}\\
\widetilde{H} & \text{for $s\in [ 1, \infty)$}.
\end{cases}
$$
Given a monotone homotopy $\widetilde{G}_s$, let
$\widetilde{J}_{s,t}$ be a family of almost complex structures in
$\widetilde{\Jj}_{[0, \, 1/2]}$ which is independent of $s$ for all
$|s|>1$ such that $\widetilde{J}_{-1,t}$ is regular for
$\widetilde{G}$ and $\widetilde{J}_{1,t}$ is regular for
$\widetilde{H}$. For each $\widetilde{x} \in
\Crit(\widetilde{\Aa}_{\widetilde{G}})$ and $\widetilde{y} \in
\Crit(\widetilde{\Aa}_{\widetilde{H}})$, let
$\wtilde{\Mm}_s(\widetilde{x},\widetilde{y}, \widetilde{G}_s,
\widetilde{J}_{s,t})$ be the moduli space of maps $\widetilde{u} \co \R
\times [0, \, 1/2] \to \widetilde{W}$ which are solutions of the
equation
\begin{equation*}
\label{s-moduli}
\partial_s \widetilde{u} + \widetilde{J}_{s,t} (\widetilde{u}) \partial_t \widetilde{u}=0,
\end{equation*}
satisfy the boundary conditions
\begin{equation*}
 \widetilde{u}(s,0) \in \Delta \text{   and   }\widetilde{u}(s,1/2) \in (\phi_{\widetilde{G}_s}^1)^{-1} (\Delta),
\end{equation*}
and have uniform limits as in \eqref{eq:moduli-limits}.
For a regular family $\widetilde{J}_{s,t}$, each moduli space
$\wtilde{\Mm}_s(\widetilde{x},\widetilde{y}, \widetilde{G}_s,
\widetilde{J}_{s,t})$ is a  smooth manifold of dimension
$\mu_{\Maslov}(\widetilde{x},\widetilde{G}) -
\mu_{\Maslov}(\widetilde{y},\widetilde{H})$. The map
$\sigma_{\widetilde{H} \widetilde{G}}$ is then defined on each
$\widetilde{x} \in \Crit^{a,b}(\widetilde{\Aa}_{\widetilde{G}})$ by
$$
\sigma _{\widetilde{H} \widetilde{G}} (\widetilde{x}) = \sum_{\scriptsize
\begin{array}{c}
  \widetilde{y} \in \Crit^{a,b}(\widetilde{\Aa}_{\widetilde{H}}),\\
  \mu_{\Maslov}(\widetilde{y},\widetilde{H}) =
  \mu_{\Maslov}(\widetilde{x},\widetilde{G}) \end{array}} \hspace{-30pt}\#
\widetilde{\Mm}_s(\widetilde{x},\widetilde{y},\widetilde{G_s},\widetilde{J}_{s,t})
 \widetilde{y},
$$
where $\# \widetilde{\Mm}_s(\widetilde{x},\widetilde{y}, \widetilde{G}_s
\widetilde{J}_{s,t})$ is the number  of elements in the
zero--dimensional moduli space $\wtilde{\Mm}_s(\widetilde{x},\widetilde{y},
\widetilde{G}_s, \widetilde{J}_{s,t})$, modulo $2$.

To prove that $\sigma_{\widetilde{H} \widetilde{G}}$ is a chain map
between the restricted Floer complexes above, one must consider the
energy of the elements $\widetilde{u} \in
\widetilde{\Mm}_s(\widetilde{x},\widetilde{y}, \widetilde{G}_s,
\widetilde{J}_{s,t})$ which is defined by the formula
$$
E(\widetilde{u}) = \int_0^{1/2} \int _{-\infty}^{+\infty} \widetilde{
\Om}(\partial_s \widetilde{u}, \widetilde{J}_{s,t}(\widetilde{u})
\partial_s \widetilde{u})\,ds\,dt.
$$
The following useful inequality is a straight forward consequence of
Stokes' theorem and the monotonicity assumption on the homotopy
$\widetilde{G}_s$.
\begin{Lemma}
For each $\widetilde{u} \in \Mm_s(\widetilde{x},\widetilde{y}, \widetilde{G}_s,
\widetilde{J}_{s,t})$ we have:
\begin{equation}\label{eq:s-energy}
E(\widetilde{u}) \leq \widetilde{\Aa}_{\widetilde{G}} (\widetilde{x}) -
\widetilde{\Aa}_{\widetilde{H}}(\widetilde{y}).
\end{equation}
\end{Lemma}
Since $E(\widetilde{u}) \geq 0$, inequality \eqref{eq:s-energy} implies
that the action must decrease under
$\sigma_{\widetilde{H}\widetilde{G}}$. This fact allows one to
verify that $\sigma_{\widetilde{H}\widetilde{G}}$ is a chain map
between the restricted Floer complexes.
The homomorphism induced  at the level of homology is also denoted
by $\sigma_{\widetilde{H}\widetilde{G}}$. It is independent of the
choices of the monotone homotopy $\widetilde{G}_s$ as well as the
family of almost complex structures $\widetilde{J}_{s,t}$.

\begin{Remark}
\label{rm:sigma} In the description of the map
$\sigma_{\widetilde{H}\widetilde{G}}$ we have tacitly assumed the
nondegeneracy of the elements of
$\Crit^{a,b}(\widetilde{\Aa}_{\widetilde{G}})$ and
$\Crit^{a,b}(\widetilde{\Aa}_{\widetilde{H}})$. If some of these
intersection points are degenerate, one simply approximates
$\widetilde{H}$ and $\widetilde{G}$ by functions $\widetilde{H}'$
and $\widetilde{G}'$ which have the nondegeneracy property and then
defines the map  between their Lagrangian Floer complexes, as above.
By the discussion at the end of Section~\ref{sec:fc}, the map induced at
the homology level is independent of the approximations
$\widetilde{H}'$ and $\widetilde{G}'$.
\end{Remark}

The following important property of monotone Floer continuation maps
is proved using \eqref{eq:s-energy} and the usual compactness and
gluing theorems from Floer theory.
\begin{Lemma}
\label{a1} Monotone Floer continuation maps  satisfy the
identities
\begin{align*}
\sigma_{\widetilde{H}\widetilde{G}} \circ
\sigma_{\widetilde{G}\widetilde{F}}& =
\sigma_{\widetilde{H}\widetilde{F}} \text{ for } \widetilde{F} \geq
\widetilde{G}
\geq \widetilde{H},\\
\sigma_{\widetilde{G}\widetilde{G}} &= \id  \text{ for every
 } \widetilde{G} \in \widetilde{\Hh}^{a,b}.
\end{align*}
\end{Lemma}

\subsection{Hamiltonian Floer homology}
Let $\Hh = C_c^{\infty}(S^1 \times W)$ be the space of smooth
time--periodic Hamiltonians on $W$ which have compact support. For
$G \in \Hh$, let $\Jj_{S^1}= \Jj_{S^1}(G)$ be the set of
$S^1$--families of $\Om$--tame almost complex structures on $W$
which are $\Om$--compatible near the support of $G$ and are equal to
$J_{gb}$ outside some compact set. In complete analogy with
Section~\ref{sec:fc}, one can associate to a generic pair $(G, J_t)
\in \Hh
\times \Jj_{S^1}$ and constants $0<a<b$,
 a restricted Hamiltonian Floer chain complex,
$$
(\CF^{a,b}(G), \partial_{J_t}).
$$
We very briefly describe this complex here and discuss how it can be
identified with a Lagrangian Floer chain complex.

Each Hamiltonian $G \in \Hh$ determines an action functional $\Aa_G
\co \Ll \to \R$ on $\Ll$, the set of smooth contractible loops in
$W$. This is given by the formula
$$
\Aa_G(x)= \int_0^1 G(x(t),t) \, dt - \int_{D^2} v^* \Omega
$$
where $v \in C^{\infty}(D^2,W)$ satisfies $v|_{\partial D^2}=x$. The
set of critical points $\Crit (\Aa_G)$ is precisely the set of
contractible $1$--periodic orbits of $G$, and each nondegenerate
critical point $x \in \Crit(\Aa_G)$ has a well-defined
Conley--Zehnder index, $\mu_{\CZ}(x,G)$. To accommodate our choice
of action functional, we normalize the  Conley--Zehnder index here
so that a local maximum of a $C^2$--small autonomous Hamiltonian on
$W^{2l}$ has Conley--Zehnder index equal to its Morse index minus
$l$.

Let $\Hh^{a,b} \subset \Hh$ be the set of all $G \in \Hh$ for which
$a$ and $b$ are not critical values of $\Aa_G$. Let
$\Crit^{a,b}(\Aa_G)$  denote the set of critical points of $\Aa_G$
with action in the interval $(a,b)$. For a generic $G \in \Hh^{a,b}$
the elements of $\Crit^{a,b}(\Aa_G)$ are nondegenerate and the
restricted Floer chain group for $G$ is defined to be the finite
dimensional, graded $\Z_2$--vector space
$$
\CF^{a,b}(G)= \bigoplus_{x \in \Crit^{a,b}(\Aa_G)} \Z_2 x.
$$
For a family $J_t\in \Jj_{S^1}$ and a pair of critical points $x,y
\in \Crit^{a,b}(\Aa_G)$, the Floer  moduli space $\Mm(x,y, G, J_t)$
consist of the maps $u \co \R \times S^1 \to M$ that satisfy the
equation
\begin{equation*}
    \partial_s u = J_t(u)(X_G(u)-\partial_t u),
\end{equation*}
and have the following limits which are uniform in $t$,
\begin{equation*}
  \lim_{s \to -\infty}u(s,t) = x(t) \qquad\text{and}\qquad
  \lim_{s \to \infty}u(s,t) = y(t).
\end{equation*}
If $J_t$ is regular, then these moduli spaces are smooth manifolds
of dimension $\mu_{\CZ}(x,G) - \mu_{\CZ}(y,G)$ and the Floer
boundary operator $\partial_{J_t} \co \CF^{a,b}(G) \to \CF^{a,b}(G)$
is defined on the generators $x \in \Crit^{a,b}(\Aa_G)$ by
$$
\partial_{J_t} (x) =
  \sum_{\scriptsize\begin{array}{c} y \in \Crit^{a,b}(\Aa_G),\\
    \mu_{\CZ}(y,G) = \mu_{\CZ}(x,G)-1 \end{array}}
  \hspace{-30pt}\# (\Mm(x,y,G,J_t)/\R) y.
$$
Here,  $\#(\Mm(x,y, G, J_t) / \R)$ is the number, modulo $2$, of
elements in the zero--dimensional space $\Mm(x,y, G, J_t)/\R$ where
$\R$ acts freely by translation in the $s$--variable.

\subsection{Translating from Hamiltonian to Lagrangian Floer theory}
\label{translation}

Given a Hamiltonian $G$ on $W$ we define a corresponding
Hamiltonian $\widetilde{G}^W $ on $\widetilde{W}$ by the formula
$$
\widetilde{G}^W(t, x_0,x_1) = G(t,x_0).
$$
One can easily check that $(\phi^1_{\widetilde{G}^W})^{-1}(\Delta)$
is equal to the graph of $\phi^1_G$ in $\widetilde{W}$. Note also
that
$$
G \geq H \Leftrightarrow \widetilde{G}^W \geq \widetilde{H}^W.
$$
Now consider the map $\Psi_{G} \co \Ll \to
\widetilde{\Ll}=\widetilde{\Ll}(\widetilde{G}^W)$ defined by
\begin{equation*}
    x(t) \mapsto \widetilde{x}(t) \eqdef \left( (\phi^t_G)^{-1}(x(t)), \phi^1_G \circ
    (\phi^{1-t}_G)^{-1}(x(1-t))\right),
\end{equation*}
where the domain of the image is restricted to $[0,\, 1/2]$.
The composition of $\Psi_G$ with $\smash{\widetilde{\Aa}_{\widetilde{G}}}$
is equal to $\Aa_G$, ie $\Aa_G(x)=
\smash{\widetilde{\Aa}_{\widetilde{G}^W}}(\widetilde{x})$.  In addition,
$\Psi_G$ identifies the sets  $\Crit(\Aa_G)$ and
$\Crit(\widetilde{\Aa}_{\widetilde{G}^W})$, and preserves
indices. More precisely, $x$ is in $\Crit(\Aa_G)$ if and only if $\widetilde{x}$
is in $\Crit(\widetilde{\Aa}_{\widetilde{G}^W})$, and  if $x$ is nondegenerate
then $\widetilde{x}$ is also nondegenerate with
$\mu_{\Maslov}(\widetilde{x},\widetilde{G}^W)=\mu_{\CZ}(x,G)$.

For a family $J_t \in \Jj_{S^1}(G)$, consider the family of
$\wtilde{\Om}$--compatible almost complex structures
$\wtilde{J}^W_{t \in [0, \, 1/2]}$  defined  by
$$
\wtilde{J}_t^W = \left((\phi^t_G)^* J_t\right) \oplus \left(-(\phi^{1-t}_G \circ
(\phi^1_G)^{-1})^* J_{1-t}\right).
$$
The map
\begin{equation*}
   u(s,t) \mapsto  \widetilde{u}(s,t) \eqdef \left( (\phi^t_G)^{-1}(u(s,t)), \phi^1_G \circ
    (\phi^{1-t}_G)^{-1}(u(s,1-t))\right),
\end{equation*}
induced by $\Psi_G$, is then a bijection
between the moduli spaces $\Mm(x,y,G,J_t)$ and $\wtilde{\Mm}(\widetilde{x},\widetilde{y},\widetilde{G}^W,\wtilde{J}^W_t)$,
for every pair $x,y \in \Crit(\Aa_G)$.
Thus, for a generic pair $(G,J_t)$, the map $\Psi_{G}$ yields an
identification between the Hamiltonian Floer complex $(\CF^{a,b}(G),
\partial_{J_t})$ and the corresponding Lagrangian Floer complex $(\CF^{a,b}
(\Delta, \widetilde{G}^W),
 \partial_{\widetilde{J}_t^W}).$
This identification preserves both the action values and indices and
yields the following isomorphism in homology
\begin{equation}\label{eq:hf=hfl}
 \HF_*^{a,b}(G) = \HF_*^{a,b}(\Delta, \widetilde{G}^W).
\end{equation}

\subsection{Morse--Bott submanifolds of periodic orbits} \label{sec:mb}

As in the Lagrangian setting, the homology groups $\HF^{a,b}(G)$
corresponding to the complex $ (\CF^{a,b}(G,J_t),
\partial_{J_t})$ are independent of the choice of the
regular family $J_t$ and are locally constant on $\Hh^{a,b}$. Again,
this allows us to define $\HF^{a,b}(G)$ for any $G \in \Hh^{a,b}$,
by approximating $G$ by a Hamiltonian whose relevant $1$--periodic
orbits are nondegenerate.

A subset $N \subset \Crit(\Aa_G)$ is said to be a \emph{Morse--Bott
manifold of periodic orbits} if the set $C_0 = \{x(0) \mid x \in
N\}$ is a closed submanifold of $W$ and $T_{x_0}C_0 = \ker (D
\phi_G^1(x_0) - Id)$ for every $x_0 \in C_0.$

It follows from \cite{bps} that a Morse--Bott manifold $ N \subset
\Ll$ of periodic orbits is mapped by $\Psi_{G}$  to a clean
intersection $\widetilde{N}$ of $\widetilde{G}^W$. The critical
submanifold also has a relative Conley--Zehnder index
$\mu_{\CZ}(N,G)$ which, by \eqref{eq:hf=hfl} and \eqref{eq:hfl=m},
is equal to $\mu_{\Maslov}(\widetilde{N},\widetilde{G}^W)$.

\subsection{Monotone Floer continuation maps for Hamiltonian
Floer homology}

For $G,H \in \Hh^{a,b}$ with $G \geq H$ one can define a monotone
chain map
$$
\sigma_{HG} \co \CF^{a,b}(G) \to \CF^{a,b}(H),
$$
in the same manner as described in Section~\ref{sec:mon-hom}. The
translation map $\Psi_G$ from  Section~\ref{translation} can again be
used to identify the maps $\sigma_{HG}$ and $\sigma_{\widetilde{H}^W
\widetilde{G}^W}$ at the chain level. In particular, we have
\begin{equation}\label{eq:switch}
 \sigma_{\widetilde{H}^W \widetilde{G}^W} \circ \Psi_{G}=
\Psi_{G} \circ \sigma_{HG}.
\end{equation}
Therefore, the maps $\sigma_{HG}$ satisfy the following version of
Lemma \ref{a1}.

\begin{Lemma}
\label{a1'} Monotone Floer continuation maps in Hamiltonian Floer
theory satisfy the identities
\begin{align*}
\sigma_{HG} \circ \sigma_{GF}& = \sigma_{\HF} \text{ for } F \geq G
\geq H,\\
\sigma_{GG} &= \id  \text{ for every
 } G \in \Hh^{a,b}.
\end{align*}
\end{Lemma}

The following additional results concerning monotone Floer
continuation maps for Hamiltonian Floer homology are well known;
see, for example, \cite{cfh,fh} and \cite[Sections~4.4--4.5]{bps}.

\begin{Lemma}
\label{a2} For constants $0< a < b < c $ and a function $G \in
\Hh^{a,b} \cap \Hh^{b,c}$ the short exact sequence of complexes
$$
0\;\ra \;\CF^{a,b}(G)\; \ra \;\CF^{a,c}(G)\; \ra \;\CF^{b,c}(G)\;
\ra\; 0,
$$
yields the exact homology triangle $\bigtriangleup_{a,b,c}(G)$ given
by
$$
\xymatrix@C=0pt{
     \HF^{a,b}(G) \ar[rr]& & \HF^{a,c}(G)\ar[dl]^{\Pi} \\
&          \HF^{b,c}(G) \ar[ul]^{\partial^*}  & }
$$
\end{Lemma}

\begin{Lemma}\label{a3}
For any $0<a<b<c$ we have the following commuting diagram
\begin{equation*}
\xymatrix{
\HF^{a,c}(G) \ar[r]^{\sigma_{HG}} &  \HF^{a,c}(H)\\
\HF^{a,b}(G)
\ar[u] 
\ar[r]^{\sigma_{HG}} & \HF^{a,b}(H)
\ar[u]. 
}
\end{equation*}
where the vertical arrows are determined by the homology triangles
$\bigtriangleup_{a,b,c}(G)$ and $\bigtriangleup_{a,b,c}(H)$.
\end{Lemma}

\begin{Lemma}
\label{a4} If $G_s$ is a monotone homotopy from $G$ to $H$ which
satisfies $G_s \in \Hh^{a,b}$ for all $s \in [-1,1]$, then
$\sigma_{HG}$ is an isomorphism.
\end{Lemma}
This last result states that the only way in which the map
$\sigma_{HG}$ can fail to be an isomorphism is if there is some $s
\in (-1,1)$ such that $G_s$ has a $1$--periodic orbit with action
equal to $a$ or $b$. Even if periodic orbits with action equal to
$a$ or $b$ appear during a monotone homotopy, it may still be
possible to show that  $\sigma_{HG}$ is nontrivial by  considering
the indices of these orbits. We now describe a useful refinement of
Lemma \ref{a4} which lies at the heart of the calculations in
\cite{cgk,fhw,gg}.

\begin{Definition}\label{def:transversal}
A monotone homotopy $G_s$ from $G$ to $H$ is said to be
\emph{transversal to} $(a,b)$ if the following conditions hold.

\begin{enumerate}
\item $G_s$ has $1$--periodic orbits with action equal to
$a$ or $b$ for only a finite set of values $\{s_j\} \subset (-1,1)$.
\item At each $s_j$, these orbits form a Morse--Bott
nondegenerate submanifold of $1$--periodic orbits, $N_j$.
\item There is an interval $I_j=(s_j-\epsilon_j, s_j+\epsilon_j)$
such that $N_j$ belongs to a smooth family of Morse--Bott
nondegenerate submanifolds  $N_j(s)$ of $G_s$ for $s \in I_j$.
\item The function
$\Aa_{G_s}(N_j(s))$ is continuous and strictly monotone on $I_j$.
\end{enumerate}

\end{Definition}
The space of such homotopies (with $\dim(N_j)=0$) is dense in the
space of all smooth monotone homotopies from $G$ to $H$.

\begin{Proposition}
\label{prop:key} Let $G_s$ be a monotone homotopy from $G$ to $H$
which is transversal to $(a,b)$. Suppose that the relative
Conley--Zehnder index of each Morse--Bott nondegenerate submanifold
$N_j$ is either strictly less than $n_0{-}\dim(N_j){-}1$ or strictly
greater than $n_0{+}1$. Then the map
$$
\sigma_{HG} \co \HF^{a,b}_{n_0}(G) \to \HF^{a,b}_{n_0}(H)
$$
is an isomorphism.
\end{Proposition}

\begin{proof}

In view of Lemma \ref{a1'}, we can clearly assume that there is only
one $s_j$. Consider then, a monotone homotopy $G_s$ from $G$ to $H$
which is transverse to $(a,b)$ such that only $G_{s_1}$ has
$1$--periodic orbits with action equal to $a$ or $b$, and these
orbits form a Morse--Bott nondegenerate submanifold $N_1$. We assume
that $\Aa_{G_{s_1}}(N_1)=b$ and the function $\Aa_{G_{s}}(N_1(s))$
is decreasing near $s_1$. The other cases can be dealt with in a
similar manner.

The set of critical values of $\Aa_{G_{s_1}}$, $\Ss(G_{s_1}) \subset
\R$, is closed and nowhere dense \cite{sc}. It follows from the
transversality assumption on $G_s$ that there is a $\delta>0$ such
that the only critical points of $\Aa_{G_{s_1}}$ with action in the
interval $[b-\delta, b]$ are those belonging to $N_1$.

We also note that the subsets  $\Ss({G_s})$ are lower
semi-continuous in $s$ in the following sense: For every open
neighborhood  $V \subset \R$ of $\Ss(G_{s'})$ there is an interval
$(s'-\epsilon , s'+ \epsilon)$ such that $\Ss_{G_s} \subset V$ for
all $s \in (s'-\epsilon , s'+ \epsilon)$ (see \cite[Section~4.4]{bps}).

It follows that we can choose an $\epsilon>0$ such that for each $s
\in [s_1-\eps, s_1 + \eps]$ the only critical points of
$\Aa_{G_{s}}$ with action in the interval $[b-\delta, b]$ are those
belonging to $N_1(s)$. By the transversality assumption, we may also
choose $\epsilon$ so that $ \Aa_{G_s}(N_1(s)) > b-\delta$ for all $s
\in [s_1-\eps, s_1 + \eps]$.

We now factor the map $\sigma_{HG}$ as
$$
\sigma_{HG} = (\sigma_{H, G_{s_1 - \epsilon}}) \circ
(\sigma_{G_{s_1-\epsilon} G_{s_1 + \epsilon}})
 \circ (\sigma_{G_{s_1+\epsilon} G}).
$$
The first and last terms are isomorphisms by Lemma \ref{a4}. It
remains to prove that the middle map is an isomorphism. To do this
we consider the commutative diagram from Lemma \ref{a3}:
\begin{equation}\label{square}
\xymatrix{ \HF^{a,b}(G_{s_1 + \epsilon})
\ar[rr]^-{\sigma_{G_{s_1-\epsilon} G_{s_1 + \epsilon}}} &
 & \HF^{a,b}(G_{s_1-\epsilon})\\
\HF^{a,b-\delta}(G_{s_1+\epsilon})
\ar[u] 
\ar[rr]^{\sigma_{G_{s_1-\epsilon} G_{s_1 + \epsilon}}} & &
\HF^{a,b-\delta}(G_{s_1-\epsilon})
\ar[u]. 
}
\end{equation}
The transversality assumption together with  our choices of $\delta$
and $\epsilon$ imply that $G_s \in \Hh^{a,b-\delta}$ for all $s \in
[s_1-\eps, s_1 + \eps]$. Hence, the bottom map is an isomorphism by
Lemma \ref{a4}.

The vertical maps belong to the exact sequences determined by Lemma
\ref{a2} for $a < b-\delta<b$. The surrounding terms in these
sequences are $\HF^{b-\delta,b}_{n_0}(G_{s_1 \pm \epsilon})$ and
$\HF^{b-\delta,b}_{n_0+1}(G_{s_1 \pm \epsilon})$. Applying Theorem
\ref{thm:poz}, we have
$$
\HF^{b-\delta,b}_{n_0}(G_{s_1 \pm \epsilon}) = H_{n_0 -
\mu_{\CZ}(N_1,G_{s_1})} (N_1)
$$
$$
\HF^{b-\delta,b}_{n_0+1}(G_{s_1 \pm \epsilon}) = H_{n_0+1 -
\mu_{\CZ}(N_1,G_{s_1})} (N_1).\leqno{\hbox{and}}
$$
The assumption that the relative index of $N_1$ is greater than
$n_0+1$ implies that these groups vanish and so the vertical arrows
in \eqref{square} are also isomorphisms.
\end{proof}

\section{The proof of Theorem~\ref{thm:main-1dyn}}

For a test Hamiltonian $H \in \Hh_{\test}(U_R)$,  Theorem
\ref{thm:main-1dyn} asserts that if $R>0$ is sufficiently small and
$\max(H)> \pi R^2$, then $H$ has a $1$--periodic orbit with action
in the interval $(\max(H), \max(H)+2\pi R^2)$. To prove this, we
will squeeze $H$ between two simple Hamiltonians  and study the
monotone Floer continuation map between their restricted
(Lagrangian) Floer homology groups.

\begin{Proposition}
\label{prop:transfer} If $R>0$ is sufficiently small, then there are
Hamiltonians  $G_+, G_- \in \Hh$ satisfying $G_+ \geq H \geq G_-$,
and constants $a,b$ satisfying
$$\max(H) < a < b < \max(H) +  2\pi R^2,$$
such that the map
$$
\sigma_{\widetilde{G}^W_-  \widetilde{G}^W_+} \co \HF^{a,b}(\Delta,
\widetilde{G}^W_+) \to \HF^{a,b}(\Delta, \widetilde{G}^W_-)
$$
is nontrivial.
\end{Proposition}

By Lemma \ref{a2} the homotopy homomorphism
$\sigma_{\widetilde{G}^W_- \widetilde{G}^W_+}$ factors as
$$
\sigma_{\widetilde{G}^W_-  \widetilde{G}^W_+} =
\sigma_{\widetilde{G}^W_-  \widetilde{H}^W} \circ
\sigma_{\widetilde{H}^W \widetilde{G}^W_+}.
$$
Hence, Proposition \ref{prop:transfer} implies that
$\HF^{a,b}(\Delta, \widetilde{H}^W)$ is nontrivial and consequently
$\CF^{a,b}(\Delta,\widetilde{H}^W)$ is nonempty. By the translation
described in Section~\ref{translation}, this, in turn, implies that
$\CF^{a,b}(H)$ is nonempty for some $a > \max(H)$. The orbit(s)
generating $\CF^{a,b}(H)$ must  be nonconstant and so Theorem
\ref{thm:main-1dyn} would follow.

To prove Proposition \ref{prop:transfer}, we construct a third model
Hamiltonian $G_0 \in \Hh$ such that $G_+  \geq G_0 \geq G_-$, and we
consider the decomposition
$$
\sigma_{\widetilde{G}^W_-  \widetilde{G}^W_+} =
\sigma_{\widetilde{G}^W_-  \widetilde{G}^W_0} \circ
\sigma_{\widetilde{G}^W_0 \widetilde{G}^W_+}.
$$
We then prove that there are constants $a$ and $b$, satisfying
$$
\max(H) <a < b < \max(H) +  2\pi R^2,
$$
and a fixed degree $n_0$, such that in this degree
$\sigma_{\widetilde{G}^W_0 \widetilde{G}^W_+}$ is a nontrivial
surjection and $\sigma_{\widetilde{G}^W_- \widetilde{G}^W_0}$ is an
isomorphism (see Propositions \ref{0-} and \ref{+0}).

\subsection{The model Hamiltonians}
We now construct the model Hamiltonians $G_+$,$G_0$, and $G_-$, and
describe their Hamiltonian flows.

\subsubsection{The function $G_+$}
\label{G-plus}

\begin{figure}[ht!]\anchor{fig:functions}
\begin{center}
{\small
\psfrag{max G_+}[][][0.8]{$\max(G_+)$} \psfrag{max
H}[][][0.8]{$\max(H)$} \psfrag{max G_-}[][][0.8]{$\max(G_-)$}
\psfrag{$c_+^1$}[][][0.8]{$c_+^1$}
\psfrag{$c_+^2$}[][][0.8]{$c_+^2$}
\psfrag{$c_-^1$}[][][0.8]{$c_-^1$}
\psfrag{$c_-^2$}[][][0.8]{$c_-^2$}
\psfrag{$d_+^1$}[][][0.8]{$d_+^1$}
\psfrag{$d_+^2$}[][][0.8]{$d_+^2$}
\psfrag{$d_-^1$}[][][0.8]{$d_-^1$}
\psfrag{$d_-^2$}[][][0.8]{$d_-^2$} \psfrag{$G_+$}[][][0.8]{$G_+$}
\psfrag{$H$}[][][0.8]{$H$} \psfrag{$G_-$}[][][0.8]{$G_-$}
\psfrag{$R$}[][][0.8]{$R$} \psfrag{$R'$}[][][0.8]{$R'$}
\psfrag{$r$}[][][0.8]{$r$}
\includegraphics[scale=0.9]{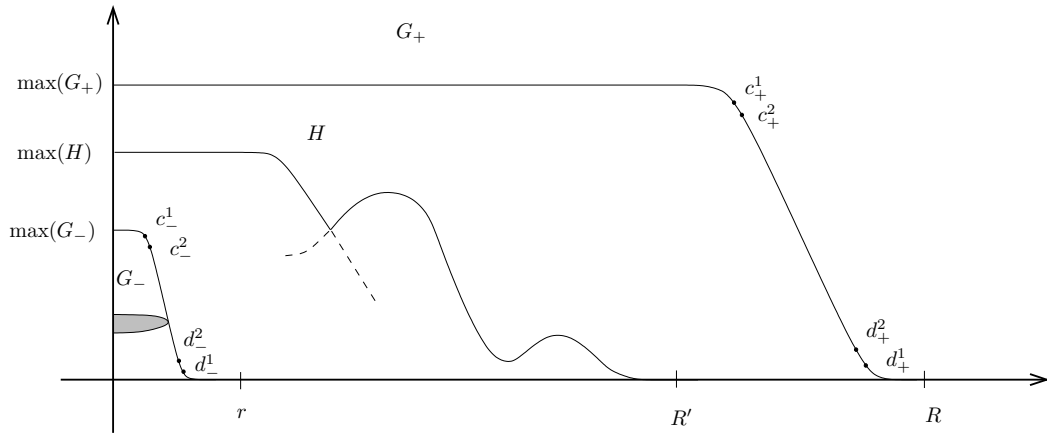}}
\caption{The functions $G_+$ and $G_-$}
\label{fig:functions}
\end{center}
\end{figure}

The function $G_+$ is constructed to approximate $H$, from above,
near its maximum set. Referring back to the notation of Section~\ref{sec:U_R}, $G_+$ is defined as a reparametrization of the
function $\|z\|^2$ on $U_R$. In particular, we set
$$
G_+(p,z) = \alpha(\|z\|^2)
$$
for a smooth nonincreasing function $\alpha$. Since $H$ is in
$\Hh_{\test}(U_R)$, there is an $R'<R$ such that $H$ vanishes for
$\|z\|>R'$. We choose $\alpha$ so that $G_+$ is constant and equal
to its maximum until $\|z\| > R'$. We then force $G_+$ to decrease
rapidly to zero at $\|z\|=R$ by choosing the slope of $\alpha$ to
decrease rapidly, become constant and finally increase to zero at
$R^2$.  The value of $\max(G_+)$ is chosen to be arbitrarily close
to, but greater than, $\max(H)$. See \figref{fig:functions}.


Since $G_+$ is a function of $\|z\|^2$, its Hamiltonian flow is just
a rescaling of the $\pi$--periodic flow of $\|z\|^2$. More
precisely, the orbits of $G_+$ on the level $\|z\|^2=c$ are all
periodic with period $\pi / \alpha'(c)$. When the slope of $\alpha$
is constant and negative we choose it not to be an integer multiple
of $\pi$. Hence, the only nonconstant $1$--periodic orbits of $G_+$
are located on two finite sequences of level sets which are labeled
in \figref{fig:functions}. One of these sequences of level sets is located in the
region where the slope decreases and the corresponding $G_+$--values
are denoted by
$$
c_+^1 > c_+^2 > \cdots
$$
The other sequence of $1$--periodic level sets is located in the
region where the slope increases. It corresponds to the
$G_+$--values
$$
d_+^1 > d_+^2 > \cdots
$$

\subsubsection{The function $G_0$}
\label{G-zero}

\begin{figure}[ht!]\anchor{fig:supports}
\begin{center}
\psfrag{M}[][][0.8]{$M$} \psfrag{pfix}[][][0.8]{$p_{\fix}$}
\psfrag{Brho}[][][0.8]{$\Bb(\rho)$} \psfrag{Br}[][][0.8]{$\Bb(r)$}
\psfrag{BR}[][][0.8]{$\Bb(R)$} \psfrag{UR}[][][0.8]{$U_R$}
\psfrag{H=maxH}[][][0.8]{$\{H= \max(H)\}$}
\psfrag{G+=c}[][][0.8]{$\{G_+=c\}$}
\psfrag{G0=c}[][][0.8]{$\{G_0=c\}$}
\psfrag{SupportG+}[][][0.8]{support of $G_+$} \psfrag{support
G0}[][][0.8]{support of $G_0$} \psfrag{G-=c}[][][0.8]{$\{G_-=c\}$}
\includegraphics{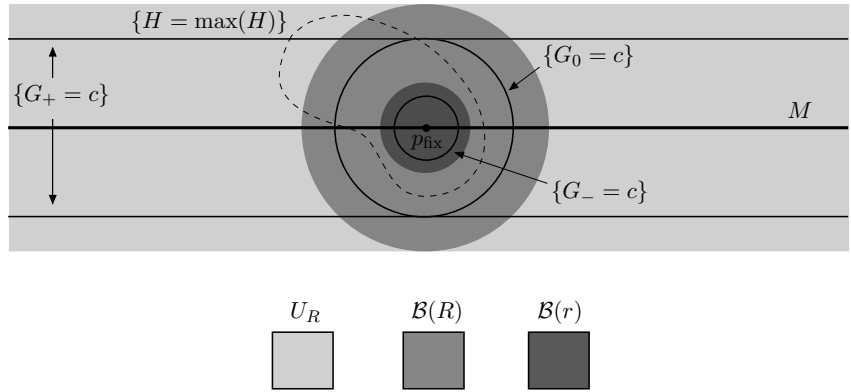}
\caption{The supports of $G_+$, $G_0$ and $G_-$}
\label{fig:supports} 
\end{center}
\end{figure}

To define $G_0$, we first describe the set on which it is supported.
Let  $p_{\fix}$  be a point on $M \subset W$. We may assume, without
loss of generality, that the interior of the set $\{ H = \max(H) \}$
intersects the zero section of $U_R$ at $(p_{\fix},0)$. This follows
from the general fact that given two points of a symplectic manifold
and a path between them, there is a Hamiltonian diffeomorphism whose
support is close to the path and whose time one flow takes one point
to the other along the path. Let $(B^{2m}(p_{\fix},\rho),
\Omega_{2m})$ be a Darboux ball in $M$ with radius $\rho$ and center
at the point $p_{\fix}$. Here, $\Om_{2m}$ denotes the standard
symplectic form on $\R^{2m}$.
The bundle $\pi \co E \to M$ is trivial over $B^{2m}(p_{\fix},\rho)$
and we may assume that the connection $\nabla$ is flat over
$B^{2m}(p_{\fix},\rho)$. Hence, there is a trivialization of $E$
over $B^{2m}(p_{\fix},\rho)$ for which the symplectic structure
$\Om$ on
 $\pi^{-1}(B^{2m}(p_{\fix},\rho)) \cap U_R$ has the form
\begin{equation}\label{Om-local}
    \Om= \Om_{2m} \oplus \Om_{2n}.
\end{equation}
Assume that $R$ is less than $\rho$.\footnote{This is the second
condition restricting the size of $R$.} We then define $G_0$ by the
formula
\begin{equation*}
G_0(p,z) = \alpha \left( \|p\|_{2m}^2 + \|z\|^2 \right),
\end{equation*}
where $\|~\|_{2m}$ is the standard norm on $B^{2m}(p_{\fix},\rho)$,
and $\alpha$ is the same function used to define $G_+$. Clearly,
$G_+ \geq G_0$ and the only place where the functions are equal and
both nonzero is on the fibre of $U_R$ over $p_{\fix}$. Moreover,
$G_0$ is supported in the Darboux ball $\Bb(R) \eqdef B^{2(m+n)}(
(p_{\fix},0),R) \subset U_R$. See \figref{fig:supports}.

With the local normal form \eqref{Om-local} for $\Om$, it is also
clear that the Hamiltonian flow of the function $\|p\|_{2m}^2 +
\|z\|^2 $ with respect to $\Om$ is $\pi$--periodic on $\Bb(R)$.
Hence, the flow of $G_0$ is again just a rescaling of a
$\pi$--periodic flow. As above, the nonconstant $1$--periodic orbits
of $G_0$ are located on two finite sequences of level sets. The
$G_0$--values of these sequences are denoted by
$$
c_0^1 > c_0^2 > \cdots
\qquad\text{and}\qquad
d_0^1 > d_0^2 > \cdots
$$

\subsubsection{The function $G_-$}
\label{G-minus} The function $G_-$ is constructed to approximate
$H$, from below, near the point $(p_{\fix},0)$ which we chose to
belong to the interior of the set $\{H = \max(H)\}$. Let $\Bb(r)$ be
a Darboux ball with radius $r < R$ and center $(p_{\fix},0)$, where
$r$ is small enough so that $\Bb(r) \subset \{H = \max(H)\}$. We
define $G_-$  by the formula
\begin{equation*} G_-(p,z) = \alpha_- ( \|p\|_{2m}^2 + \|z\|^2).
\end{equation*}
The function $\alpha_-$ has the same general behavior as $\alpha$.
It is chosen so that $G_-$ is equal to its maximum in a small
neighborhood of $(p_{\fix},0)$ and then decreases to zero within
$\Bb(r)$.  When the slope of $\alpha_-$ is constant and negative we
choose it not to be an integer multiple of $\pi$. We also choose
$\alpha_-$ so that $\max(G_-)$ is arbitrarily close to, but smaller
than, $\max(H)$.
 For these choices we have  $H \geq G_-$ (see \figref{fig:functions}).

Again, the nonconstant $1$--periodic orbits of $G_-$ occur on the
level sets corresponding to two finite sequences of $G_-$--values,
as pictured in \figref{fig:functions}. We label the two sequences  by
$$
c_-^1 > c_-^2 > \cdots
\qquad\text{and}\qquad
d_-^1 > d_-^2 > \cdots
$$

\subsubsection{Morse--Bott nondegeneracy}

The functions $\alpha$ and $\alpha_-$ used to define $G_+$, $G_0$,
and $G_-$ are chosen so that their second derivatives are not zero
when their slope is nonconstant. This implies that all the
nonconstant $1$--periodic level sets are Morse--Bott nondegenerate
as defined in Section~\ref{sec:mb} (see \cite[Lemma~5.3.2]{bps}).

\subsubsection{Actions and indices}

Table \ref{action-index} summarizes the necessary information
concerning the periodic orbits of our model Hamiltonians. For
convenience, we will identify a level set $\{G =c \}$ comprised of
$1$--periodic orbits with its corresponding value $c$.  The
$1$--periodic orbits on such a level, all have the same action. This
common action value is given in the third column of Table
\ref{action-index} up to a term which is denoted by ellipses
``$\cdots$'' and can be made arbitrarily small.

\def\strutt{\vrule width 0pt height 14pt depth 6pt}
\begin{table}[h!]
  \centering
\begin{tabular}{|c|c|c|}
  \hline
 \strutt Level set, clean int.  & Relative index & Action \\
  \hline
  \strutt$c_+^k,\,\widetilde{c}_+^k$ & $(2k-1)n-m+1$ & $\max(G_+) + k \pi R^2 + \cdots$ \\
  \strutt$d_+^k,\,\widetilde{d}_+^k$ & $(2k-1)n-m$ & $k \pi R^2 + \cdots$ \\
  \hline
  \strutt$c_0^k,\,\widetilde{c}_0^k$ & $(2k-1)(m+n)+1$ & $\max(G_0) +k \pi R^2 + \cdots$ \\
  \strutt$d_0^k,\,\widetilde{d}_0^k$ & $(2k-1)(m+n)$ & $k \pi R^2 + \cdots$ \\
  \hline
  \strutt$c_-^k,\,\widetilde{c}_-^k$ & $(2k-1)(m+n)+1$ & $\max(G_-) + k \pi r^2 + \cdots$ \\
  \strutt$d_-^k,\,\widetilde{d}_-^k$ & $(2k-1)(m+n)$ & $k \pi r^2 + \cdots$ \\
  \hline
\end{tabular}
\caption{Actions and indices}\label{action-index}
\end{table}

Each of the level sets of $1$--periodic orbits is nondegenerate in
the Morse--Bott sense and so has a relative Conley--Zehnder index.
These indices are listed in the second column of
Table~\ref{action-index}. For a discussion of the calculation of the relative
indices of the level sets of type $c_+^k$ and $d_+^k$, the reader is
referred to \cite{cgk,gg}. By construction, the levels of type $c_0^k$,
$d_0^k$, $c_-^k$ and $d_-^k$ can be treated as if they are contained in
$\R^{2(m+n)}$ and the relative indices can be easily derived from the
calculations in \cite{fhw}.

To prove Proposition \ref{prop:transfer} we also need to utilize the
corresponding objects in the Lagrangian setting which are obtained
using the translation described in Section~\ref{translation}. For the
translations of the model Hamiltonians we still have
$$
\widetilde{G}^W_+ \geq \widetilde{G}^W_0 \geq \widetilde{G}^W_-.
$$
Since the $1$--periodic level sets of the model Hamiltonians are all
Morse--Bott nondegenerate, they each get mapped to a clean
intersection in $\widetilde{W}$. We denote the  clean intersection
corresponding to the level $c$ by $\widetilde{c}$, and note that $c$ and
$\widetilde{c}$ are  diffeomorphic and have
 the same action value and relative index.
For example, the map $\Psi_{G_+}$  takes $c_+^1$ to a clean
intersection $\widetilde{c}_+^1$ of $\widetilde{G}^W_+$, such that every
$1$--periodic orbit $x(t) \subset c_+^1$ gets mapped to the
intersection point $(x(0),x(0)) \in \Delta \cap
(\phi^1_{\widetilde{G}^W_+})^{-1}(\Delta)$.

\subsection{The constants $a$ and $b$}

We choose the constants $a<b$ so that
\begin{align}\label{ab}
\max(H) &< a < \max(G_-)+\pi r^2 < \max(G_+) +\pi R^2 \notag\\
        &< b < \max(G_+) + 2\pi R^2.
\end{align}
Since $\max(H)>\pi R^2$, the following periodic level sets do not
have actions in the interval $(a,b)$ and may be ignored: $d^1_+$,
$d^1_0$, $d^1_-$, and the levels $c^k_+$ and $c^k_0$ for $k \geq 2$.

\subsection[The homomorphisms sigma and sigma]{The homomorphisms $\sigma_{\widetilde{G}^W_-  \widetilde{G}^W_0}$ and
$\sigma_{\widetilde{G}^W_0  \widetilde{G}^W_+}$}

\begin{Proposition}
\label{0-} For $a<b$ as in \eqref{ab} and $n_0 = m+n+1$, we have
$$\HF^{a,b}_{n_0}(\Delta, \widetilde{G}^W_0) =
\HF^{a,b}_{n_0}(\Delta, \widetilde{G}^W_-) = \Z_2,$$ and the
homomorphism
$$
\sigma_{\widetilde{G}^W_-  \widetilde{G}^W_0} \co
\HF^{a,b}_{n_0}(\Delta, \widetilde{G}^W_0) \to
\HF^{a,b}_{n_0}(\Delta, \widetilde{G}^W_-)
$$
is an isomorphism.
\end{Proposition}

\begin{Proposition}
\label{+0} For $a<b$ as in \eqref{ab} and $n_0=m+n+1$, the
homomorphism
$$
\sigma_{\widetilde{G}^W_0  \widetilde{G}^W_+} \co
\HF^{a,b}_{n_0}(\Delta, \widetilde{G}^W_+) \to
\HF^{a,b}_{n_0}(\Delta, \widetilde{G}^W_0)
$$
is surjective.
\end{Proposition}

Together, Propositions \ref{0-} and \ref{+0} imply that the map
$$
\sigma_{\widetilde{G}^W_- \widetilde{G}^W_+} =
\sigma_{\widetilde{G}^W_- \widetilde{G}^W_0} \circ
\sigma_{\widetilde{G}^W_0 \widetilde{G}^W_+} \co
\HF^{a,b}_{n_0}(\Delta, \widetilde{G}^W_+) \to
\HF^{a,b}_{n_0}(\Delta, \widetilde{G}^W_-)
$$
is nontrivial and Proposition \ref{prop:transfer} follows.\qed

\subsection{Proof of Proposition~\ref{0-}}

By the identifications \eqref{eq:hf=hfl} and \eqref{eq:switch} from
Section~\ref{translation}, it suffices for us to prove Proposition
\ref{0-} for the corresponding Hamiltonian Floer homology groups.
That is, we may prove that
$$
\HF^{a,b}_{n_0}(G_0) = \HF^{a,b}_{n_0}(G_-) = \Z_2,
$$
and the map
$$
\sigma_{G_- G_0} \co \HF^{a,b}_{n_0}(G_0) \to \HF^{a,b}_{n_0}(G_-)
$$
is an isomorphism.

First we consider the group  $\HF_{n_0}^{a,b}(G_0) =
\HF_{n_0}^{a,b}(\Delta, \widetilde{G}_0)$. By our choice of the
constants $a<b$ from \eqref{ab}, we know that
$\Crit^{a,b}(\widetilde{\Aa}_{\widetilde{G}^W_0})$ consists of the
points in $\widetilde{c}^1_0$ and possibly one other  clean intersection
of type  $\widetilde{d}^k_0$ for some $k \geq 2$. We will show that the group
$HF_{n_0}^{a,b}(\Delta, \widetilde{G}_0)$ is determined  solely by
$\widetilde{c}^1_0$. In particular, we show that
$HF_{n_0}^{a,b}(\Delta, \widetilde{G}_0)$ corresponds, via Theorem
\ref{thm:poz}, to $H_0(\widetilde{c}^1_0; \Z_2)$. The group will be shown
 not to depend on the intersections of type $\widetilde{d}^k_0$ with $k \geq 2$,
 because their relative Maslov indices are too large.

We first note that we may choose the function $G_0$ such that the
action of $\widetilde{c}^1_0$ is distinct from the actions of the $\widetilde{d}^k_0$,
for $k \geq 2$. We may then choose constants $a'$ and $b'$
with $a < a' < \widetilde{\Aa}_{\widetilde{G}_0}(\widetilde{c}^1_0) < b'
<b$, such that
\begin{equation*}
\Crit^{a',b'}(\widetilde{\Aa}_{\widetilde{G}^W_0})=\widetilde{c}^1_0.
\end{equation*}
Applying  Theorem \ref{thm:poz} to the clean intersection
$\widetilde{c}^1_0$, we get
$$\HF^{a',b'}_{n_0} (\Delta, \widetilde{G}_0) = H_{0}(\widetilde{c}_0^1;\Z_2)
= H_{0}(S^{2(m+n)-1}; \Z_2)
= \Z_2.$$
To prove that $\HF_{n_0}^{a,b}(G_0) = \Z_2$, it remains for us to
show that
\begin{equation*}
\label{ends} \HF^{a',b'}_{n_0} (\Delta,
\widetilde{G}_0)=\HF^{a,b}_{n_0} (\Delta, \widetilde{G}_0).
\end{equation*}
Using  the exact triangle from Lemma \ref{a2} (twice),  we see that
it suffices to show that the groups $\HF^{a,a'}_* (\Delta,
\widetilde{G}_0)$ and $\HF^{b,b'}_{*+1} (\Delta, \widetilde{G}_0)$
are trivial for $*=n_0$ and $*= n_0-1$.

As described above, for the action window determined by $a$ and $b$,
 the only way in which the groups $\HF^{a,a'}_* (\Delta, \widetilde{G}_0)$ and $\HF^{b,b'}_{*+1} (\Delta, \widetilde{G}_0)$
 can be  nontrivial is if they are generated by a clean intersection of type $\widetilde{d}^k_0$ with $k \geq 2$.  In fact, only one of the
 intervals $(a,\,a')$ and $(b\,,b')$ can include the action of some $\widetilde{d}^k_0$.
So, by Theorem \ref{thm:poz}, we either have
$$
\HF^{a,a'}_*(\Delta, \widetilde{G}_0) = H_{* -
\mu_{\Maslov}(\widetilde{d}_0^k, \widetilde{G}_0)}(S^{2(m+n)-1} ; \Z_2)
$$
$$
 \HF^{b,b'}_{*+1}(\Delta, \widetilde{G}_0) = H_{* +1- \mu_{\Maslov}(\widetilde{d}_0^k, \widetilde{G})}(S^{2(m+n)-1} ; \Z_2).\leqno{\hbox{or}}
$$
Now, the relative index of each $\widetilde{d}_0^k$ with $k \geq 2$
satisfies $$\mu_{\Maslov}(\widetilde{d}_0^k, \widetilde{G}_0) \geq
3n_0-3 \geq 2n_0.$$ Thus, for
$*= n_0$ and $n_0-1$,  these groups are indeed trivial  as required.

The fact that $HF^{a,b}_{n_0}(G_-)=  \HF^{a,b}_{n_0}(\Delta,
\widetilde{G}_-)= \Z_2$,  can be established using a similar
argument. The only difference in this case is that the groups
$\HF^{a,a'}_*(\Delta, \widetilde{G}_-)$ and  $\HF^{b,b'}_{*+1}(\Delta,
\widetilde{G}_-)$ may both be nontrivial. However they are still generated
by clean intersections of type $\widetilde{d}^k_-$ for $k \geq 2$, and so for
$*=n_0$ and $n_0-1$ these groups still vanish.

\begin{figure}[ht!]\anchor{fig:homotopy}
\begin{center}
\psfrag{R}[][][0.8]{$R$} \psfrag{r}[][][0.8]{$r$}
\psfrag{maxG-}[][][0.6]{$\max(G_-)$}
\psfrag{maxG0}[][][0.6]{$\max(G_0)$} \psfrag{G-}[][][0.8]{$G_-$}
\psfrag{G0}[][][0.8]{$G_0$}
\includegraphics[scale=0.9]{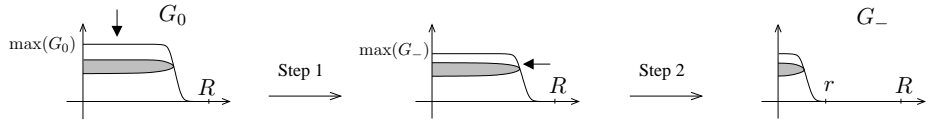}
\caption{The homotopy in two steps}
\label{fig:homotopy} 
\end{center}
\end{figure}

To prove that the map  $\sigma_{G_- G_0}$ is an isomorphism in
degree $n_0$ we consider the homomorphism induced by the monotone
homotopy shown in \figref{fig:homotopy}. This homotopy is transversal to the
interval $(a,b)$ as per Definition \ref{def:transversal}, and the
Morse--Bott nondegenerate submanifolds of periodic orbits, $N_i$,
are spheres of dimension $2(m+n)-1$.

For simplicity, the total monotone homotopy is broken into two
steps. This allows us to factor the total Floer continuation
homomorphism $\sigma_{G_- G_0}$ as the composition of two
homomorphisms. We prove that each of these factors is an isomorphism
in degree $n_0$.


During the first homotopy, the maximum of the function $G_0$ is
decreased to coincide with $\max(G_-)$. This is achieved by
decreasing the absolute value of the constant slope. Since
$\max(G_0)$ and $\max(G_-)$ are arbitrarily close, it is clear that
the level sets of $1$--periodic orbits change only slightly during
this homotopy and no new periodic orbits with period one and action
equal to $a$ or $b$ are created. By Lemma \ref{a4} the homomorphism
induced by this step of the homotopy is an isomorphism in
Hamiltonian Floer homology restricted to actions in $(a,b)$.

During the second homotopy the radius of the ``spike" is decreased
but its slope remains constant (see \figref{fig:homotopy}). More precisely, the
Hamiltonians  which occur in this homotopy are of the form
$G_s=\alpha_s(\|p\|_{2m}^2+\|z\|^2)$, and the only change which
occurs is that the size of the interval on which  $\alpha_s $ take
its maximum value decreases. Hence, each of the Hamiltonians  $G_s$
has its $1$--periodic orbits located on two finite sequences of
$1$--periodic level sets, $\{c^k(s)\}$ and $\{d^k(s)\}$. The
relative indices of the these level sets are
$$
\mu_{\CZ}(c^k(s),G_s)=(2k-1)(m+n)+1
$$
and
$$
\mu_{\CZ}(d^k(s),G_s)=(2k-1)(m+n).
$$
We note that for $k \geq 2$ all of the relative indices are greater
than $n_0+1 = m+n+2$.

The actions $\Aa_{G_s}(c^k(s))$  and $\Aa_{G_s}(d^k(s))$ are smooth
decreasing functions of $s$. We now determine which of the levels
can attain the action values $a$ or $b$. The functions
$\Aa_{G_s}(c^k(s))$ are always greater than $a$ and can only attain
the value $b$ when $k \geq 2$, since, for $k=1$, the action starts
at $\Aa_{G_{-1}}(c^1(-1)) \approx \Aa_{G_0}(c^1_0)=\max(G_+)+ \pi
R^2 +\dots<b$. The functions $\Aa_{G_s}(d^k(s))$ can only take the
values $a$ or $b$ for $k \geq 2$ since, for $k=1$, the action starts
at $\Aa_{G_{-1}}(d^1(-1)) \approx \Aa_{G_{0}}(d^1_0)=\pi R^2+ \dots
<a$. As noted above, the relative indices for $k \geq 2$  are all
greater than $m+n+2$, so it follows from Proposition \ref{prop:key}
that the monotone Floer continuation map induced by the second
monotone homotopy is an isomorphism in degree $n_0 = m+n+1$. This
completes the proof of Proposition \ref{0-}.\qed

\subsection{Proof of Proposition~\ref{+0}}

For $\star = 0  \text{ and }+$, we know that the set
$\Crit^{a,b}(\smash{\widetilde{\Aa}_{\widetilde{G}_{\star}}})$ consists of
$\widetilde{c}^1_{\star}$ and possibly one other clean intersection of type
$\widetilde{d}^k_{\star}$ with $k \geq 2$. In considering the map
$\sigma_{\widetilde{G}^W_0 \widetilde{G}^W_+}$ acting on
$\HF_{n_0}^{a,b}(\delta, \widetilde{G}_+)$,  the clean intersections
of type $\widetilde{d}^k_{\star}$ may again be avoided as they were in
 the proof of Proposition \ref{0-}. More precisely, we may assume that the actions of $\widetilde{c}^1_{\star}$ are
 distinct from the actions of the $\widetilde{d}^k_{\star}$ for $k \geq 2$, and choose  constants
$a'$ and $b'$ with $a < a' <
\widetilde{\Aa}_{\widetilde{G}_{\star}}(\widetilde{c}^1_{\star}) < b'
<b$, such that
\begin{equation*}
\Crit^{a',b'}(\widetilde{\Aa}_{\widetilde{G}^W_{\star}})=\widetilde{c}^1_{\star}.
\end{equation*}
The same arguments used to prove Proposition \ref{0-}, now imply that the
groups $\HF^{a,a'}_* (\Delta,
\widetilde{G}_{\star})$ and $\HF^{b,b'}_{*+1} (\Delta, \widetilde{G}_{\star})$
are trivial for $*=n_0$ and $*= n_0-1$. It then follows from repeated use of Lemma \ref{a2} and Lemma
\ref{a3}, that proving Proposition \ref{+0} it equivalent to showing  that
the map
$$
\sigma_{\widetilde{G}^W_0  \widetilde{G}^W_+} \co
\HF^{a',b'}_{n_0}(\Delta, \widetilde{G}^W_+) \to
\HF^{a',b'}_{n_0}(\Delta, \widetilde{G}^W_0)
$$
is surjective. To prove this, we will need to study the map
$\sigma_{\widetilde{G}^W_0 \widetilde{G}^W_+}$ at the chain level.
As described in Remark \ref{rm:sigma}, this requires us to consider
approximations $\widetilde{G}^{W'}_+ $ of $\widetilde{G}^W_+$ and
${\widetilde{G}^{W'}_0}$ of $\widetilde{G}^W_0$, such that the
$1$--periodic orbits of $\widetilde{G}^{W'}_+$ and
$\widetilde{G}^{W'}_0$ with actions in $(a',b')$ are nondegenerate.
\begin{Proposition}
\label{x-star} The approximations $\widetilde{G}^{W'}_+$ and
$\widetilde{G}^{W'}_0$ can be constructed so that there is a point
$\widetilde{x}_*$ in $\widetilde{W}$ with the following properties:
\begin{enumerate}
  \item The point $\widetilde{x}_*$ lies in the intersection
    $\Crit_{n_0}^{a',b'}(\smash{\widetilde{\Aa}_{\widetilde{G}^{W'}_+}})
    \cap \Crit_{n_0}^{a',b'}(\smash{\widetilde{\Aa}_{\widetilde{G}^{W'}_0}})$.
  \item The class $[\widetilde{x}_* + \widetilde{w}]$ is nontrivial in
    $\HF^{a',b'}_{n_0}(\Delta,\widetilde{G}^{W'}_+)$ for some chain
    $\widetilde{w}$ in $\CF_{n_0}^{a',b'}(\Delta,\widetilde{G}^{W'}_+)$
    with $\smash{\Aa_{\widetilde{G}^{W'}_+}}(\widetilde{w}) <
    \smash{\Aa_{\widetilde{G}^{W'}_+}}(\widetilde{x}_*)$.
  \item The class $[\widetilde{x}_*]$ is nontrivial and generates
    $\HF^{a',b'}_{n_0}(\Delta,\widetilde{G}^{W'}_0)$.
  \item At the chain level, $\sigma_{\widetilde{G}^{W'}_0
    \widetilde{G}^{W'}_+}(\widetilde{x}_* + \widetilde{w}) =
    \widetilde{x}_*$.
\end{enumerate}
\end{Proposition}

At the level of homology, Proposition \ref{x-star} implies that
$$
\sigma_{\widetilde{G}^{W'}_0 \widetilde{G}^{W'}_+}([\widetilde{x}_* +
\widetilde{w}]) = [\widetilde{x}_*] \neq 0.
$$
In particular, the map
$$
\sigma_{\widetilde{G}^{W'}_0 \widetilde{G}^{W'}_+} \co
\HF^{a',b'}_{n_0}(\Delta, \widetilde{G}^{W'}_+) \to
\HF^{a',b'}_{n_0}(\Delta, \widetilde{G}^{W'}_0)
$$
is surjective. By Remark \ref{rm:sigma}, the same conclusion
holds for $\sigma_{\widetilde{G}^W_0 \widetilde{G}^W_+}$, and
by the discussion above this implies Proposition \ref{+0}.\qed

\subsection{Proof of Proposition~\ref{x-star}}

For our choice of $a'$ and $b'$, we have
\begin{equation*}\label{eq:c1plus}
\widetilde{c}^1_+ = \Crit^{a',b'}(\smash{\widetilde{\Aa}_{\widetilde{G}^W_+}})
\qquad\text{and}\qquad \widetilde{c}^1_0 =
\Crit^{a',b'}(\smash{\widetilde{\Aa}_{\widetilde{G}^W_0}}).
\end{equation*}
Let $h_+ \co \widetilde{c}^1_+ \to \R$ be a Morse function and let $g^+$
be a metric on $\widetilde{c}^1_+$ such that the corresponding Morse
complex is well-defined. Similarly choose a Morse function and
metric $(h_0,g_0)$ on $\widetilde{c}^1_0$. Theorem \ref{thm:poz} then
yields, for sufficiently small $\delta>0$, the following strong
equivalences of chain complexes:
\begin{equation}\label{plus}
(C(h_+),\, \partial_{g_+}) \equiv
(\CF^{a',b'}(\Delta, \widetilde{G}^W_+ + \delta
\widetilde{H}_{h_+}), \,
\partial_{J_{g_+,t}} )
\end{equation}
and
\begin{equation}\label{0}
(C(h_0),\, \partial_{g_0}) \equiv
(\CF^{a',b'}(\Delta, \widetilde{G}^W_0 + \delta
\widetilde{H}_{h_0}), \,
\partial_{J_{g_0,t}}).
\end{equation}
We recall here the important point that this strong equivalence
implies that the critical points of the Morse functions are exactly
the points of intersection which generate the corresponding
Lagrangian Floer complexes. The equivalences \eqref{plus} and
\eqref{0} will be used to prove Proposition \ref{x-star} at the
level of Morse homology. In particular, we will set
$\widetilde{G}^{W'}_+ = \widetilde{G}^W_+ + \delta
\widetilde{H}_{h_+}$ and $\widetilde{G}^{W'}_0 = \widetilde{G}^W_0 +
\delta \widetilde{H}_{h_0}$ for special choices of the Morse
functions $h_+$ and $h_0$.

The definition of these functions involves the following simple
model. Let
$$
\height_k \co S^{k} \subset \R^{k+1} \to \R
$$ denote the height function
on the standard $k$--dimensional unit sphere given by the
restriction of the function
$$
(x_1,\dots, x_{k+1}) \mapsto x_1.
$$
For a (linear) subsphere of $S^k$ defined by $l$ linear conditions
of the form $x_j=0$ where $j \neq 1$, the restriction of $\height_k$
clearly agrees with the function $\height_{k-l}$ and the critical
points of the two functions coincide.

Recall that $\widetilde{c}_+^1$ is an $S^{2n-1}$--bundle over $M$. To
define $h_+$ we first choose a Morse function $f_M \co M \to \R$
which, for simplicity, we assume is self-indexing and has only one
critical point of index $2m$ at $p_{\fix} \in M$. In other words, if
$\Crit(f_M)=\{p_j\}$, then
\begin{equation*}\label{f-plus}
f_M(p_j) = \mu_{\operatorname{Morse}}(p_j, h_+), \text{   and    } (f_M)^{-1}(2m)=
p_{2m} = p_{\fix}.
\end{equation*}
The lift of $f_M$ to $\widetilde{c}_+^1$ is a Morse--Bott function whose
critical submanifolds correspond to the fibres of $\widetilde{c}_+^1$
over the points $p_j$. To define $h_+$ we perturb the lift of $f_M$
near each of its critical fibres as follows
\begin{equation*}\label{eq:h-plus}
h_+(p,z)= f_M(p) + \epsilon \sum_{p_j \in \Crit(f_M)} \eta_j(p)
\height_{2n-1}(z).
\end{equation*}
Here, each $\eta_j$ is a smooth bump function on $M$ which has
support near $p_j$ and attains its maximum value, one, in a
neighborhood of $p_j$. For sufficiently small $\epsilon>0$, $h_+$ is
a Morse function whose critical points lie in the fibres over the
points $p_j \in \Crit(f_M)$ where they coincide with the critical
points of the function $\height_{2n-1}$ on the fibre. In other words,
each critical point $p_j$ of $f_M$ gives rise to exactly two
critical points of the perturbation $h_+$, which we denote by
$x_{j,+}^{\xtop}$ and $x_{j,+}^{\xbottom}$. The Morse index of
$x_j^{\xtop}$ is given by
$$
\mu_{\operatorname{Morse}}(x_{j,+}^{\xtop}) =\mu_{\operatorname{Morse}}(p_j) + 2n-1
$$
and similarly
$$
\mu_{\operatorname{Morse}}(x_{j,+}^{\xbottom}) =\mu_{\operatorname{Morse}}(p_j).
$$
Note, that the function $h_+$ restricts to the fibre of
$\widetilde{c}_+^1$ over $p_j$ as
\begin{equation}\label{eq:h-plus-restriction}
f_M (p_j)+\epsilon \cdot \height_{2n-1}.
\end{equation}
Now, the level set $\widetilde{c}_0^1$ is diffeomorphic to
$S^{2(m+n)-1}$, and, by the definitions of $G_+$ and $G_0$, it is
clear that $\widetilde{c}_+^1 \cap \widetilde{c}_0^1$ is equal to the fibre
of $\widetilde{c}_+^1$ over $p_{\fix} \in M$. We define the Morse
function $h_0 \co \widetilde{c}_0^1 \approx S^{2(m+n)-1} \to \R$ by the
equation
\begin{equation*}\label{h-0}
    h_0 = f_M(p_{\fix}) + \epsilon \cdot \height_{2(m+n)-1}.
\end{equation*}
As described above, the restriction of $h_0$ to $\widetilde{c}_+^1 \cap
\widetilde{c}_0^1$ agrees with the restriction of $h_+$ (see
\eqref{eq:h-plus-restriction}) and the two critical points of $h_0$,
which we call $x_0^{\xtop}$ and $x_0^{\xbottom}$, coincide with
$x_{2m,+}^{\xtop}$ and $x_{2m,+}^{\xbottom}$, respectively.

\subsubsection{Proof of part (1)} Set
$$
\widetilde{x}_* = x_{2m,+}^{\xbottom} = x_0^{\xbottom}.
$$
Under the identifications \eqref{plus} and \eqref{0}, it follows
that the point $\widetilde{x}_* \in \wtilde{W}$ belongs to both
$\Crit^{a',b'}(\smash{\widetilde{\Aa}_{\widetilde{G}^{W'}_+}})$ and
$\Crit^{a',b'}(\smash{\widetilde{\Aa}_{\widetilde{G}^{W'}_0}})$.

The critical point $x_{2m,+}^{\xbottom}$ has Morse index $2m$ and the
clean intersection $\widetilde{c}_+^1$ has relative Maslov index
$n-m+1$. It follows from Theorem \ref{thm:poz} that $\widetilde{x}_*$
has Maslov index $(2m)+(n-m+1) = m+n+1$ as an element of
$\Crit^{a',b'}(\smash{\widetilde{\Aa}_{\widetilde{G}^{W'}_+}})$. Similarly,
$\widetilde{x}_*$ has Maslov index $m+n+1$ as an element of
$\Crit^{a',b'}(\smash{\widetilde{\Aa}_{\widetilde{G}^{W'}_0}})$.\qed

\subsubsection{Proof of part (2)}

Let $g_+$ be a metric on $\widetilde{c}_+^1$ for which the Morse chain
complex $(C(h_+),\partial_{g_+})$ is well-defined. This complex has
a natural filtration coming from the function $f_M$ on $M$. More
precisely, let $C_{i,j}(h_+)$ be the $\Z_2$--vector space generated
by the critical points $x$ of $h_+$ such that
$\mu_{\operatorname{Morse}}(x,h_+)=i+j$ and $x$ lies in the fibre
over some $p_k \in \Crit(f_M)$  where
$\mu_{\operatorname{Morse}}(p_k,f_M)=i$. Then
\begin{equation*}
C_k=\bigoplus_{i+j=k} C_{i,j}
\end{equation*}
admits the filtration
\begin{equation*}
    \Ff_0 C_k \subseteq \Ff_1 C_k \dots\subseteq \Ff_{k} C_k,
\end{equation*}
where
\begin{equation*}
 \Ff_{j} C_k = \bigoplus_{i \leq j} C_{i,k-i}.
\end{equation*}
Since we are only considering coefficients in the field $\Z_2$, the
corresponding spectral sequence converges to $H_*(\widetilde{c}_+^1,
\Z_2)$. It is easy to see that the $E^1$--term is given by
$$
E^1_{i,j}=C_i(f_M) \otimes H_j(S^{2n-1}, \Z_2) = C_{i,j}.
$$
\begin{Lemma}
\label{lem:e2} The $E^2$--term of the spectral sequence is
$$
E^2_{i,j}=H_i(M,H_j(S^{2n-1},\Z_2)).
$$
\end{Lemma}
A proof of this well-known result has been included in Appendix B.

To prove Part (2) of Proposition \ref{x-star} we first note that
$E^1_{2m,0} = \Z_2$ and is represented by $\widetilde{x}_* =
x_{2m,+}^{\xbottom}$. Let $\beta_k$ denote the $k^{th}$ $\Z_2$--Betti
number of $M$. It follows from Lemma \ref{lem:e2} that
\begin{equation}\label{e2m}
E^2_{2m} = E^2_{2m,0} \oplus E^2_{2(m-n)+1,2n-1}= \Z_2  \oplus
\Z_2^{\beta_{2(m-n)+1}}.
\end{equation}
In addition, the assumption that the unit normal bundle of $M$ is
homologically trivial in degree $2m$ implies that
\begin{equation}\label{h2m}
H_{2m}(\widetilde{c}_+^1, \Z_2) = \Z_2  \oplus \Z_2^{\beta_{2(m-n)+1}}.
\end{equation}
Since the spectral sequence converges to $H_*(\widetilde{c}_+^1, \Z_2)$,
equations \eqref{e2m} and \eqref{h2m} imply that $E^k_{2m,0}$ is
isomorphic to $\Z_2$ and is generated by $[\widetilde{x}_*]$ for all $k
\geq 0$. Moreover, $\widetilde{x}_*$ must appear nontrivially in a
representative of a nonzero class in $H_{2m}(\widetilde{c}_+^1, \Z_2)$.
We denote this class by $[\widetilde{x}_* + \widetilde{w}]$

By the strong equivalence of complexes from \eqref{plus} and our
choice of $h^+$ it is clear that $\widetilde{x}_*$ has the largest
action in
$\Crit_{n_0}^{a',b'}(\widetilde{\Aa}_{\widetilde{G}^{W'}_+})$. This
completes the proof of Part (2).\qed

\subsubsection{Proof of part (3)}

The point $\widetilde{x}_* = x^{\xbottom}_0$ is the unique global minimum
of $h_0$. Since
$$
H_*(C(h_0), \p_{g_0}) = H_*(\widetilde{c}_0^1, \Z_2)=H_*(S^{2(m+n)-1},
\Z_2),
$$
the class $[\widetilde{x}_*]$ must represent the nontrivial class
$H_0(C(h_0), \p_{g_0})$. Part (3) then follows immediately from the
strong equivalence of complexes from \eqref{0}.\qed

\subsubsection{Proof of part (4)}

Let $\widetilde{v}_+ \in \Crit^{a',b'}(\widetilde{G}^{W'}_+)$ and
$\widetilde{v}_0\in \Crit^{a',b'}(\widetilde{G}^{W'}_0)$. For a monotone
homotopy $\widetilde{G}^{W'}_s$ from $\widetilde{G}^{W'}_+$ to
$\widetilde{G}^{W'}_0$ and a generic family of almost complex
structures $J_{s,t}$, consider an element $\widetilde{u} \in
\smash{\wtilde{\Mm}_s(\widetilde{v}_+,\widetilde{v}_0,\widetilde{G}^{W'}_s,\wtilde{J}_{s,t})}$.
By inequality \eqref{eq:s-energy}, we have
$$
\widetilde{\Aa}_{\widetilde{G}^{W'}_0}(\widetilde{v}_0) \leq
\widetilde{\Aa}_{\widetilde{G}^{W'}_+}(\widetilde{v}_+)
$$
with equality only when $\widetilde{u}(s,t)$ does not depend on $s$.

Setting $\widetilde{v}_+ = \widetilde{x}_* = \widetilde{v}_0$ and noting that
$\smash{\widetilde{\Aa}_{\widetilde{G}^{W'}_+}}(\widetilde{x}_*)=
\smash{\widetilde{\Aa}_{\widetilde{G}^{W'}_0}}(\widetilde{x}_*)$ it follows that
the moduli space
$\smash{\wtilde{\Mm}_s(\widetilde{x}_*,\widetilde{x}_*,\widetilde{G}^{W'}_s,\wtilde{J}_{s,t})}$
contains only the constant map.

We also know that $\widetilde{x}_*$ has the smallest
$(\smash{\widetilde{\Aa}_{\widetilde{G}^{W'}_0}})$--value in
$\Crit^{a',b'}(\smash{\widetilde{\Aa}_{\widetilde{G}^{W'}_0}})$. Since any
$\widetilde{v}_+$ which appears in $\widetilde{w}$ with nonzero coefficient
satisfies
$$\widetilde{\Aa}_{\widetilde{G}^{W'}_+}(\widetilde{v}_+) <
\widetilde{\Aa}_{\widetilde{G}^{W'}_0}(\widetilde{x}_*),$$
we have
$\smash{\wtilde{\Mm}_s(\widetilde{v}_+,\widetilde{x}_*,\widetilde{G}^{W'}_s,
\wtilde{J}_{s,t})} = \emptyset$.

By the definition of the map $\sigma_{\widetilde{G}^{W'}_0
\widetilde{G}^{W'}_+}$ we get $\sigma_{\widetilde{G}^{W'}_0
\widetilde{G}^{W'}_+} (\widetilde{x}_* + \widetilde{w}) = \widetilde{x}_*$, as
desired.\qed

\section*{Appendices}\small
\addcontentsline{toc}{section}{Appendices}
\appendix

\section{Area bounds for plane curves with positive curvature}
\label{app:a}

We begin by recalling the statement of Proposition \ref{prop:area}.

\begin{Proposition}\label{prop:area-here}
Let $\xi \co [0,T] \to \R^2$ be a closed planar curve with constant
speed $v$, rotation number $k$, and positive curvature $K(t)/v$ such
that $0 < \underline{K} \leq K(t)$. The Euclidean area enclosed by
$\xi$, $A(\xi)$, satisfies
$$
0 \leq A(\xi) \leq k 4(v/\underline{K})^2.
$$
\end{Proposition}

\begin{proof}

A closed piecewise linear curve is called a \emph{box curve} if each
segment is either vertical or horizontal. We can assume that as one
traverses the curve the segments alternate between vertical and
horizontal. A box curve is said to be \emph{positive} if one always
makes a {\em left turn} between segments as one follows the curve in
a counter-clockwise manner.  Note that any positive box curve is
completely determined (up to translation) by two finite sequences of
positive numbers $\{a_j\}_{j=1}^{2k}$ and $\{b_j\}_{j=1}^{2k}$ such
that
$$
\sum_{j=1}^{2k} (-1)^{j+1}a_j =0 = \sum_{j=1}^{2k} (-1)^{j+1}b_j.
$$
Here the $a_j$ are the lengths of the vertical segments and the
$b_j$ are the lengths of the horizontal segments and we use the
convention that the first segment corresponds to $a_1$. The rotation
number of the curve is clearly equal to $k$. (See \figref{fig:boxcurve}.)


\begin{figure}[ht!]\anchor{fig:boxcurve}
\begin{center}
\psfrag{a1}[][][.8]{$a_1$}
\psfrag{a2}[][][.8]{$a_2$}
\psfrag{a3}[][][.8]{$a_3$}
\psfrag{a4}[][][.8]{$a_4$}
\psfrag{a5}[][][.8]{$a_5$}
\psfrag{a6}[][][.8]{$a_6$}
\psfrag{b1}[][][.8]{$b_1$}
\psfrag{b2}[][][.8]{$b_2$}
\psfrag{b3}[][][.8]{$b_3$}
\psfrag{b4}[][][.8]{$b_4$}
\psfrag{b5}[][][.8]{$b_5$}
\psfrag{b6}[][][.8]{$b_6$}
\includegraphics{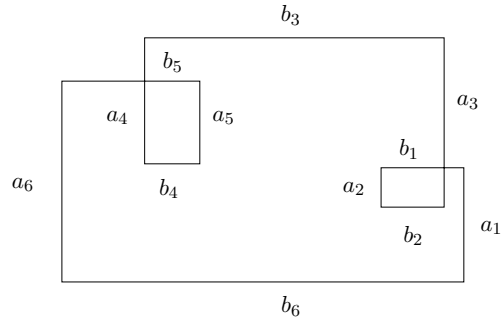}
\caption{A positive box curve}
\label{fig:boxcurve}
\end{center}
\end{figure}


To prove Proposition \ref{prop:area} we will associate to $\xi(t)$ a
positive box curve $\widehat{\xi}(t)$ such that
$$
A(\xi) \leq A(\widehat{\xi}) \leq k4 (v/\underline{K})^2.
$$
Note that
$$
\dot{\xi}(t)=\left(v \cos\theta(t), v \sin\theta(t)\right)
$$
where $\theta(t)$ is the solution of the initial value problem
$$
\dot{\theta} =K(t), \,\,\,\,\, \theta(0)=0.
$$
We define the value $t_j$ as the unique solution of the equation
$$
\theta(t_j) = j \pi.
$$
Similarly, we let $\tau_j$ be the solution of
$$
\theta(\tau_j) = {\textstyle\frac{\pi}{2}} + j \pi.
$$
For $j=1, \dots, 2k$ we then set
$$ a_j = y(t_j)-y(t_{j-1})\,\,\,\,\,  \text{ and
}\,\,\,\,\, b_j = x(t_j)-x(t_{j-1}).
$$
These sequences determine the box curve $\widehat{\xi}$  illustrated in
\figref{fig:boxarea}. Clearly,
$$
A(\widehat{\xi}) \geq A(\xi).
$$
Moreover, for all $j=1, \dots, 2k$ we have
$$
a_j,\, b_j \leq 2v/\underline{K} .
$$


\begin{figure}[ht!]\anchor{fig:boxarea}
\begin{center}
\psfrag{xi}{$\xi$}
\psfrag{xi-hat}{$\widehat{\xi}$}
\includegraphics{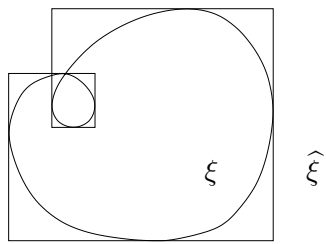}
\caption{The curves $\xi$ and $\widehat{\xi}$}
\label{fig:boxarea}
\end{center}
\end{figure}

\begin{Lemma}
\label{hat}
Let $\widehat{\xi}$ be a positive box curve determined by
the sequences $\{a_j\}_{j=1}^{2k}$ and $\{b_j\}_{j=1}^{2k}$. Suppose
that
$$
a_j,\, b_j \leq \overline{c}
$$
for all $j=1, \dots, 2k$. Then
$$
0 < A(\widehat{\xi}) \leq k\overline{c}^2.
$$
\end{Lemma}

\proof
For convenience we set
$$
\widetilde{a}_i = \sum_{j=1}^{i} (-1)^{j+1}a_j.
$$
The area of $\widehat{\xi}$ is then given by the formula
$$
A(\widehat{\xi})= b_1\widetilde{a}_1 -b_2 \widetilde{a}_2 + \,\, \dots \,\,+
b_{2k-1}\widetilde{a}_{2k-1}.
$$
Consider the positive box curve $\overline{\xi}$ determined by the
sequences $\{\overline{a}_j\}_{j=1}^{2k}$ and $\{\overline{b}_j\}_{j=1}^{2k}$
where $\overline{a}_j =a_j$ and $\overline{b}_j=\overline{c}.$ We then have
\begin{eqnarray*}
  A(\overline{\xi}) &=&  \overline{c}\widetilde{a}_1 -\overline{c}\widetilde{a}_2 + \,\, \dots \,\,+
\overline{c}\widetilde{a}_{2k-1}\\
  {} &=& \left(b_1\widetilde{a}_1 - \,\, \dots \,\,+
b_{2k-1}\widetilde{a}_{2k-1}\right) + \left((c-b_1)\widetilde{a}_1 - \,\,
\dots \,\,+ (c-b_{2k-1})\widetilde{a}_{2k-1}\right) \\
  {} &=& A(\widehat{\xi}) + \left((c-b_1)\widetilde{a}_1
- \,\, \dots \,\,+ (c-b_{2k-1})\widetilde{a}_{2k-1}\right).
\end{eqnarray*}
The second term is the area of the positive box corresponding to the
sequences $\{a_j\}_{j=1}^{2k}$ and $\{\overline{c}-b_j\}_{j=1}^{2k}$.
Hence this term is nonnegative and we have
$$
A(\overline{\xi}) \geq A(\widehat{\xi}).
$$
Finally,
$$A(\overline{\xi}) =  \overline{c}(\widetilde{a}_1-\widetilde{a}_2 +
  \cdots + \widetilde{a}_{2k-1})
  = \overline{c}(a_1 + a_3 +  \cdots + a_{2k-1}) \leq k
  \overline{c}^2.\qed$$
Applying Lemma \ref{hat} to the positive box curve $\widehat{\xi}$
defined by $\xi$, we obtain Proposition \ref{prop:area-here}.
\end{proof}

\section{Spectral sequences for the Morse homology of fibre bundles}
\label{app:b}

Let $pr \co P \to B$ be a fibre bundle with a closed base $B$ and
with closed fibres diffeomorphic to a manifold $F$. Let $h_B \co B
\to \R$ be a self-indexing Morse function with critical points
$\{b_j\}$, ie $h_B(b_j)= \mu_{\operatorname{Morse}}(b_j, h_B)$.
The lift of $h_B$ via $pr$ is a Morse--Bott function on $P$ whose
critical submanifolds are the fibres of $P$ over the points $b_j$.
We can construct from $h_B$ a useful Morse function on $P$ by
perturbing $pr^*h_B$ close to its critical fibres as follows. Let
$\Vv_j$ be a small open neighborhood of $b_j$ in $B$ over which the
bundle $P$ is trivial and let $\eta_j$ be a smooth bump function on
$B$ which is supported in $\Vv_j$ and attains its maximum in a
neighborhood of $b_j$. Define the function $h \co P \to \R$ by
\begin{equation*}\label{h}
h(b,z)= h_B(b) + \epsilon \sum_{b_j \in \Crit(h_B)} \eta_j(b)
f(z),
\end{equation*}
where $f \co F \to \R$ is chosen to be Morse. For sufficiently small
$\epsilon>0$, $h$ is Morse and all its critical points lie in the
fibres $pr^{-1}(b_j) \approx F$ where they coincide with the
critical points of the function $f$.

Fix a metric $g$ on $P$ for which the Morse chain complex
$(C(h),\partial_g)$ is well-defined. That is, we assume that the
stable and unstable manifolds of the gradient vector field
$$
V_{h,g} = \grad(h,g)
$$
intersect transversally. The chain group $C(h)$ is the
$\Z_2$--vector space generated by the critical points of $h$. The
boundary operator is defined using the moduli spaces,
$\textsl{m}(x,y)$, of integral curves of $-V_{h,g}(x)$ going from $x
\in \Crit(h)$ to $y \in \Crit(h)$. By our assumption on $V_{h,g}$,
each $\textsl{m}(x,y)$ is a smooth manifold of dimension
$\mu_{\operatorname{Morse}}(x,h) - \mu_{\operatorname{Morse}}(y,h)$
on which $\R$ acts freely by translation. The map $\p_g$ is then
defined by
$$
\p_g(x) = \sum_{\mu_{\operatorname{Morse}}(y,h)= \mu_{\operatorname{Morse}}(x,h)-1} \#
(\textsl{m} (x,y)/ \R) \cdot y,
$$
where $ \# (\textsl{m} (x,y)/ \R)$ is the number of elements in
$\textsl{m} (x,y)/\R$ modulo $2$.

Our special choice of the function $h$ yields a natural filtration
for the chain complex $(C(h),\partial_g)$. More precisely, let
$C_{i,j}(h)$ be the $\Z_2$--vector space generated by the critical
points $x$ of $h$ that satisfy
$$
\mu_{\operatorname{Morse}}(x,h)=i+j\,\,\,\, \text{   and
}\,\,\,\,\mu_{\operatorname{Morse}}(pr(x), h_B)=i.
$$
Then each graded component
\begin{equation*}
C_k(h)=\bigoplus_{i+j=k} C_{i,j}(h)
\end{equation*}
admits the filtration
\begin{equation*}
    \Ff_0 C_k(h) \subseteq \Ff_1 C_k(h) \dots\subseteq \Ff_{k} C_k(h),
\end{equation*}
where
\begin{equation*}
 \Ff_{j} C_k(h) = \bigoplus_{i \leq j} C_{i,k-i}(h).
\end{equation*}
The assumption that $h_B$ is self-indexing implies the following
useful fact: For $x \in C_{i,j}(h)$  and $y \in C_{i',j'}(h)$ the
moduli space $\textsl{m}(x,y)$ is empty  whenever  $i'>i$. This
follows from the fact that $h(x)<h(y)$. It implies that the boundary
map respects the filtration and we get a spectral sequence
$E^k_{i,j}$.

Let $\p^k \co C_{i,j}(h) \to C_{i-k,j+k-1}(h)$ be the map which is
defined on the generators $x \in C_{i,j}(h)$ by
$$
\p^k(x) =
\sum_{\scriptsize\begin{array}{c}
  \mu_{\operatorname{Morse}}(y,h)=i+j-1,\\
  \mu_{\operatorname{Morse}}(pr(y),h_B)=i-k
  \end{array}} \hspace{-30pt}\#(\textsl{m} (x,y)/\R ) \cdot y.
$$
The Morse boundary operator $\p_g$ decomposes as
$$
\p_g = \sum_{k\geq0} \p^k.
$$
Moreover, for each $k \geq 0$ we have $\p^k \circ \p^k =0$ and
$$
E_{i,j}^{k+1} = H_*(E_{i,j}^{k}, \p^k).
$$
Since we are using coefficients in the field $\Z_2$ the
spectral sequence converges to $H_*(P, \Z_2)$. That is,
$E_{i,j}^k = E_{i,j}^{k+1}=  E_{i,j}^{\infty}$ for
all sufficiently large $k$, and
$$
H_l(P, \Z_2)=\bigoplus_{i+j=l}E_{i,j}^{\infty}.
$$
It is easy to see that the $E^1$--terms are given by
$$
E^1_{i,j} = C_i(h_B) \otimes H_j(F,\Z_2).
$$
The point of this appendix is to provide a proof of the following
well-known result.
\begin{Proposition}
\label{prop:e2} The $E^2$--term of the spectral sequence is
$$
E^2_{i,j}=H_i(B,H_j(F, \Z_2)).
$$
\end{Proposition}

Our proof of Proposition \ref{prop:e2} utilizes the generalization
of Morse homology to gradient-like vector fields. This argument is
motivated by the recent work of Hutchings \cite{hutch}.

Recall that a vector field $V$ is a \emph{gradient-like} vector field
for a Morse function $h \co N^k \to \R$ if
\begin{itemize}
    \item $dh(V) > 0$ away from the critical points of $h$,
    \item near each critical point $x$ of $h$, $V$ has the
form of a negative gradient vector field $-V_{h,g}$ for some metric
$g$.
\end{itemize}
If $V$ is gradient-like for $h$ and the stable and unstable
manifolds of its flow intersect transversely, then $V$ is said to be
Morse--Smale. Given such a vector field, one can construct a Morse
complex $(C(h), \p_V)$ in the usual way.

The resulting homology $H_*(C(h),\p_V)$ does not depend on the
choice of $V$ and is equal to $H_*(N)$. To see this consider two
Morse--Smale gradient-like vector fields for $h$, $V_1$ and $V_2$.
For $s \in \R$ let $V_s$ be a family of vector fields such that $V_s
=V_1$ for all $s \leq -1$, $V_s =V_2$ for all $s \geq 1$, and
$dh(V_s)>0$ away from $\Crit(h)$ for all $s \in \R$ (such families
always exist). As in Floer theory, this determines a continuation
chain map
$$
\sigma_{V_2,V_1} \co (C(h), \p_{V_1}) \to (C(h), \p_{V_2}).
$$
The map is defined by counting
the solutions $u \co \R \to N$ of the equation
$$
\dot{u}(s) = -V_s(u(s)).
$$
The usual arguments imply that at the homology level
$\sigma_{V_3,V_2 }\circ\sigma_{V_2,V_1 }= \sigma_{V_3,V_1 }$ and
$\sigma_{V,V }= Id$. Hence $H_*(C(h),\p_V)$ does not depend on
$V$. Letting $V = \grad(h,g)$ for a generic metric $g$ on $N$ it
follows that
$$
H_*(C(h),\p_V) = H_*(N).
$$

\subsection{Proof of Proposition~\ref{prop:e2}}

Let $V_1 = V_{h,g}$. The idea of the proof of Proposition
\ref{prop:e2} is to first construct another gradient-like vector
field $V_2$ for $h$ which admits its own spectral sequence
$\Ee_{i,j}^k$ such that $\Ee_{i,j}^2=H_i(B,H_j(F, \Z_2))$. Then we
show that the continuation map $\sigma_{V_1,V_2}$ induces an
isomorphism from $E_{i,j}^2$ to $\Ee_{i,j}^2$.

To define $V_2$ we follow \cite{hutch}. Let $h_B \co B \to \R$ and
$f \co F \to \R$ be the functions used to define $h$. Let $g_B$ be a
metric on $B$ for which the gradient vector field $V_B =
\grad(h_B,g_B)$ is Morse--Smale. Let $g_F$ be a fibrewise metric on
the total space $P$ and let $V_F$ be the fibrewise gradient field of
$f$ with respect to $g_F$ . Fixing a connection on $P$ we set
\begin{equation*}\label{v2}
    V_2 = V_F + \Hor(V_B),
\end{equation*}
where $\Hor(V_B)$ is the horizontal lift of $V_B$. Clearly, $V_2$ is
gradient-like for $h$. For a generic choice of the fibrewise metric
$g_F$, $V_2$ is also Morse--Smale \cite{hutch}. In addition, we may
assume that the pairs $(f, g_F(b_j))$ are Morse--Smale for every
critical point $b_j$ of $h_B$.

The boundary map of the corresponding Morse complex $(C(h),
\p_{V_2})$ also respects the filtration
\begin{equation*}
    \Ff_0 C_k(h) \subseteq \Ff_1 C_k(h) \dots\subseteq \Ff_{k} C_k(h),
\end{equation*}
and yields a second spectral sequence $\Ee_{i,j}^k$.
\begin{Lemma}
$\Ee_{i,j}^2=H_i(B,H_j(F,\Z_2))$
\end{Lemma}

\begin{proof}
The following identification is obvious
$$
\Ee^0_{i,j}= C_{i,j}(h)=C_i(h_B) \otimes C_j(f).
$$
The boundary operator $\p_V$ decomposes as
$$
\p_V = \sum_{k\geq0} \p_V^k,
$$
where each $\p^k \co C_{i,j}(h) \to C_{i-k,j+k-1}(h)$ is defined as
above. Hence,
$$
\Ee_{i,j}^{k+1} = H_*(\Ee_{i,j}^{k}, \p^k),
$$
and we have
$$
\Ee^1_{i,j}= C_i(h_B) \otimes H_j(F,\Z_2).
$$
It suffices for us to prove that
\begin{equation}\label{last}
\p^1_V(b \otimes [\alpha]) = \p_{g_B}(b) \otimes [\alpha],
\end{equation}
for each critical point $b \in \Crit_i(h_B)$ and each class
$[\alpha] \in H_j(F, \Z_2)$.

By definition
\begin{equation*}
\p^1_V(b \otimes [\alpha]) = \sum_{c \in \Crit_{i-1}(h_B)}\qua
\sum_{[\beta] \in H_j(F, \Z_2)}
\# \left(\textsl{m}( b \otimes \alpha, c \otimes \beta)/\R \right)
c \otimes [\beta].
\end{equation*}
Consider a trajectory  $\gamma \in \textsl{m}( b \otimes \alpha, c
\otimes \beta)$. By the definition of $V_2$ the projection of
$\gamma$ to $B$ is a gradient trajectory $\gamma_B$ of $h_B$ with
respect to $g_B$. Let us choose a trivialization of $P$ over
$\gamma_B$ and note that this yields a homotopy of functions $h_s
=h|_{pr^{-1}(\gamma_B(s))}$ on $F$ and metrics $g_s=
g_F(\gamma_B(s))$ on $F$.

In our trivialization the set of trajectories of $V_2$ which project
to $\gamma_B$ are then distinguished by their fibre components. It
is easy to see that these fibre components are exactly the
trajectories counted in the Morse continuation map from $(C(f),
\p_{g_F(b)})$ to $(C(f), \p_{g_F(c)})$ that is determined by the
homotopies $h_s$ and $g_s$. Since the continuation map induces an
isomorphism in homology we have $[\alpha] = [\beta]$.

Summing over the gradient trajectories of $h_B$ from $b$ to $c$ we
obtain \eqref{last}.
\end{proof}

\begin{Lemma}
The map $\sigma_{V_1,V_2}$ is a morphism of filtered complexes and
hence it induces an isomorphism from $E_{i,j}^k$ to $\Ee_{i,j}^k$
for all $k \geq 0.$
\end{Lemma}

\begin{proof}
By construction $\sigma_{V_1,V_2}$  is a chain map. To prove that it
is a morphism of filtered complexes it remains to show that
$\sigma_{V_1,V_2}(\Ff_i C_*(h)) \subset \Ff_i C_*(h).$ Let $x \in
\Ff_i C_*(h)$ and suppose that $y$ appears in the image of $x$ under
$\sigma_{V_1,V_2}$. It follows from the definition of the map
$\sigma_{V_1,V_2}$  that $h(y) \leq h(x)$. Since $h$ is a
perturbation of the self-indexing function $h_B$ this implies that
$y \in \Ff_i C_*(h)$.
\end{proof}

\normalsize

\bibliographystyle{gtart}

\end{document}

{fig:functions}=1
{fig:supports}=2
{fig:homotopy}=3
{fig:boxcurve}=4
{fig:boxarea}=5

%% file: gtspec.tex
\def\ifplaintex{\expandafter\ifx\csname documentclass\endcsname\relax}

\def\ifplaintex{\expandafter\ifx\csname documentclass\endcsname\relax}

\ifplaintex 
\else
\fi

\expandafter\ifx\csname epsfbox\endcsname\relax\input epsf\fi

\def\gt{{\mathsurround=0pt\it $\cal G\mskip-2mu$eometry \&\ 
$\cal T\!\!$opology}}        

\def\gtp{{\mathsurround=0pt\it $\cal G\mskip-2mu$eometry \&\ 
$\cal T\!\!$opology $\cal P\!$ublications}}  

\def\lognumber#1{\def\thelognumber{#1}}
\def\volumenumber#1{\def\thevolumenumber{#1}}
\def\papernumber#1{\def\thepapernumber{#1}}
\def\volumeyear#1{\def\thevolumeyear{#1}}

\def\pagenumbers#1#2{\def\startpage{#1}\def\finishpage{#2}}
\def\published#1{\def\publishdate{#1}}
\def\proposed#1{\def\theproposer{#1}}
\def\seconded#1{\def\theseconders{#1}}
\def\received#1{\def\receiveddate{#1}}
\def\revised#1{\def\reviseddate{#1}}
\def\accepted#1{\def\accepteddate{#1}}
\def\asciititle#1{\def\theasciititle{#1}}
\def\covertitle#1{\def\thecovertitle{#1}}

\def\asciiaddress#1{\def\theasciiaddress{#1}}

\def\asciiurl#1{\def\theasciiurl{#1}}
\long\def\asciiabstract#1{\long\def\theasciiabstract{#1}}
\def\asciikeywords#1{\def\theasciikeywords{#1}}

\let\\\par\let\thelognumber\relax
\let\thevolumenumber\relax\let\thepapernumber\relax
\let\thevolumeyear\relax\let\thesamplenumber\relax\let\startpage\relax
\let\finishpage\relax\let\publishdate\relax\let\receiveddate\relax
\let\reviseddate\relax\let\accepteddate\relax\let\theasciititle\relax
\let\thecovertitle\relax\let\theasciiauthors\relax\let\theasciiaddress\relax
\let\theasciiabstract\relax\let\theasciikeywords\relax
\let\theasciiemail\relax\let\theshortauthors\relax\let\theshorttitle\relax
\let\theasciiurl\relax

\long\def\maketitlep{   

\count0=\startpage

\gt\hfill      
\hbox to 77pt{\vbox to 0pt{\vglue -15pt\epsfbox{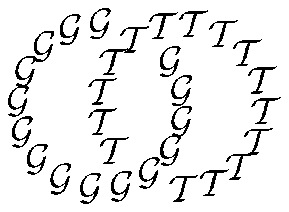}\vss}\hss}
\break
{\small\ifx\thesamplenumber\relax 
Volume \else Sample
\fi\thevolumenumber\ (\thevolumeyear)
\startpage--\finishpage\nl
Published: \publishdate}
\vglue 0.5truein plus 0.4fil minus 0.1truein

{\parskip=0pt\leftskip 0pt plus 1fil\def\\{\par\smallskip}{\ifplaintex\large
\else\Large\fi\bf\thetitle}\par\medskip}   

\vglue 0pt plus 0.1fil 

{\parskip=0pt\leftskip 0pt plus 1fil\def\\{\par}{\sc\theauthors}
\par\medskip}

\vglue 0pt plus 0.1fil 

{\small\parskip=0pt\let\newline\\
{\leftskip 0pt plus 1fil\def\\{\par}{\sl\theaddress}\par}
\expandafter\ifx\theemail\relax    
\relax\else\vglue 5pt plus 0.02fil minus 2pt\def\\{\stdspace{\rm 
and}\stdspace} 
\cl{Email:\stdspace\tt\theemail}\fi
\ifx\theurl\relax                  
\relax\else\vglue 5pt plus 0.02fil minus 2pt\def\\{\stdspace{\rm 
and}\stdspace}
\cl{URL:\stdspace\tt\theurl}\fi\par}

\vglue 7pt plus 0.3fil minus 3pt

{\bf Abstract}
\vglue 5pt plus 0.1fil minus 2pt

\theabstract

\vglue 7pt plus 0.3fil minus 3pt

{\bf AMS Classification numbers}\quad Primary:\quad \theprimaryclass

Secondary:\quad \thesecondaryclass

\vglue 5pt plus 0.3fil minus 2pt

{\bf Keywords:}\quad \thekeywords

\vglue 10pt plus 0.5fil minus 5pt

{\small  Proposed: \theproposer\hfill Received: \receiveddate\nl
Seconded: \theseconders\hfill 
\ifx\reviseddate\relax                         
Accepted: \accepteddate                        
\else
Revised: \reviseddate                          
\fi}
\eject
}       

\font\phead=cmsl9 scaled 950
\font=cmsl9 scaled 1050
\font\pnum=cmbx10 scaled 913
\font\lnum=cmbx10 
\font\pfoot=cmsl9 scaled 950
\font=cmsl9 scaled 1050
\ifplaintex
\headline{\vbox to 0pt{\vskip -4.5mm\line{\small\phead\ifnum
\count0=\startpage ISSN 1364-0380 (on line)
1465-3060 (printed) \hfill {\pnum\folio}\else\ifodd\count0\def\\{ }%
\ifx\theshorttitle\relax\thetitle\else\theshorttitle\fi\hfill{\pnum\folio}
\else\def\\{ and }{\pnum\folio}\hfill\ifx\theshortauthors\relax\theauthors
\else\theshortauthors\fi\fi\fi}\vss}}
\footline{\vbox to 0pt{\vglue 0mm\line{\small\pfoot\ifnum\count0=\startpage
\copyright\ \gtp\hfill\else
\gt, Volume \thevolumenumber\ (\thevolumeyear)\hfill\fi}\vss
}}
\else
\makeatletter
\def\@oddhead{{\small\lhead\ifnum\count0=\startpage ISSN 1364-0380 (on line)
1465-3060 (printed) \hfill {\lnum\number\count0}\else\ifodd\count0
\def\\{ }\ifx\theshorttitle\relax \thetitle \else\theshorttitle\fi\hfill
{\lnum\number\count0}\else\def\\{ and }{\lnum\number\count0}
\hfill\ifx\theshortauthors\relax 
\theauthors\else\theshortauthors\fi\fi\fi}}\def\@evenhead{\@oddhead}
\def\@oddfoot{\small\lfoot\ifnum\count0=\startpage\copyright\ \gtp\hfill\else
\gt, Volume \thevolumenumber\ (\thevolumeyear)\hfill\fi}
\def\@evenfoot{\@oddfoot}
\makeatother
\fi

\newwrite\gtoutfile
\long\gdef\makeheadfile{  
{\def\\{, }\def\s{ }
\immediate\openout\gtoutfile head.xxx
\immediate\write\gtoutfile{Proxy-for: \ifx\theasciiauthors\relax
\theauthors\else\theasciiauthors\fi\s<\ifx\theasciiemail\relax\theemail\else\theasciiemail\fi>}
\immediate\write\gtoutfile{\noexpand\\}
\immediate\write\gtoutfile{Authors: \ifx\theasciiauthors\relax
\theauthors\else\theasciiauthors\fi}
{\def\\{ }\immediate\write\gtoutfile{Title: \ifx\theasciititle\relax
\thetitle\else\theasciititle\fi}}
\immediate\write\gtoutfile{Subj-class: GT or SG or MG etc}
\immediate\write\gtoutfile{MSC-class: \theprimaryclass\ifx\thesecondaryclass\relax\else, \thesecondaryclass\fi}
\immediate\write\gtoutfile{Journal-ref: Geom. Topol. \thevolumenumber
(\thevolumeyear) \startpage-\finishpage}
\immediate\write\gtoutfile{Comments: Published by Geometry and Topology at}
\immediate\write\gtoutfile{\s\s http://www.maths.warwick.ac.uk/gt/GTVol\thevolumenumber/paper\thepapernumber.abs.html}
\immediate\write\gtoutfile{\noexpand\\}
\immediate\write\gtoutfile{}
\ifx\theasciiabstract\relax
\immediate\write\gtoutfile{\theabstract}\else
\immediate\write\gtoutfile{\theasciiabstract}\fi
\immediate\write\gtoutfile{}
\immediate\write\gtoutfile{\noexpand\\}
\immediate\write\gtoutfile{}
\immediate\closeout\gtoutfile}}  

\def\maketitlepage{\maketitlep\makeheadfile}
\let\maketitle\maketitlepage